\documentclass[11pt]{article}
\usepackage{a4wide}

\newcommand{\ba}{\begin{array}}
\newcommand{\ea}{\end{array}}
\newcommand{\be}{\begin{equation}}
\newcommand{\ee}{\end{equation}}
\newcommand{\la}{\label}
\newcommand{\bea}{\begin{eqnarray}}
\newcommand{\eea}{\end{eqnarray}}
\newcommand{\ch}{\choose}
\renewcommand{\a}{\alpha}
\renewcommand{\b}{\beta}
\renewcommand{\c}{\gamma}
\newcommand{\G}{\Gamma}
\renewcommand{\L}{L_n^{(\a)}(x)}
\newcommand{\Ln}{L_{n+1}^{(\a)}(x)}
\newcommand{\KL}{L_n^{\a,M}(x)}
\newcommand{\gL}{L_n^{\a,M,N}(x)}
\renewcommand{\l}{\left}
\renewcommand{\r}{\right}
\newcommand{\set}[1]{\left\{#1\right\}_{n=0}^{\infty}}
\newcommand{\hyp}[5]{\mbox{}_{#1}F_{#2}
\left(\left.{#3 \atop #4}\right|#5\right)}
\newcommand{\hyphyp}[5]{\mbox{}_{#1}\Phi_{#2}
\left(\left.{#3 \atop #4}\right|#5\right)}
\newcommand{\ndots}{n=0,1,2,\ldots}
\newcommand{\n}{\nonumber}
\newcommand{\nn}{\nonumber \\}
\newcommand{\ds}{\displaystyle}
\newcounter{stelling}
\newcommand{\st}[1]{\par\vspace{0.5cm}\refstepcounter{stelling}
{\bf Theorem \thestelling.} {\sl #1}\par\vspace{0.5cm}}
\newcommand{\lem}[1]{\par\vspace{0.5cm}\refstepcounter{stelling}
{\bf Lemma \thestelling.} {\sl #1}\par\vspace{0.5cm}}
\newcommand{\cor}[1]{\par\vspace{0.5cm}\refstepcounter{stelling}
{\bf Corollary \thestelling.} {\sl #1}\par\vspace{0.5cm}}

\begin{document}

\title{On differential equations for Sobolev-type Laguerre polynomials}
\author{J.~Koekoek \and R.~Koekoek \and H.~Bavinck}
\date{ }
\maketitle

\begin{abstract}
The Sobolev-type Laguerre polynomials $\set{\gL}$ are orthogonal with
respect to the inner product
$$<f,g>\;=\frac{1}{\G(\a+1)}\int_0^{\infty}x^{\a}e^{-x}f(x)g(x)dx+Mf(0)g(0)+
Nf'(0)g'(0),$$
where $\a>-1$, $M\ge 0$ and $N\ge 0$.

In 1990 the first and second author showed that in the case $M>0$ and $N=0$
the polynomials are eigenfunctions of a unique differential operator of the
form
$$M\sum_{i=1}^{\infty}a_i(x)D^i+xD^2+(\a+1-x)D,$$
where $\l\{a_i(x)\r\}_{i=1}^{\infty}$ are independent of $n$. This
differential operator is of order $2\a+4$ if $\a$ is a nonnegative integer
and of infinite order otherwise.

In this paper we construct all differential equations of the form
\bea & &M\sum_{i=0}^{\infty}a_i(x)y^{(i)}(x)+
N\sum_{i=0}^{\infty}b_i(x)y^{(i)}(x)+{}\nn
& &\hspace{1cm}{}+MN\sum_{i=0}^{\infty}c_i(x)y^{(i)}(x)+
xy''(x)+(\a+1-x)y'(x)+ny(x)=0,\n\eea
where the coefficients $\l\{a_i(x)\r\}_{i=1}^{\infty}$,
$\l\{b_i(x)\r\}_{i=1}^{\infty}$ and $\l\{c_i(x)\r\}_{i=1}^{\infty}$ are
independent of $n$ and the coefficients $a_0(x)$, $b_0(x)$ and $c_0(x)$ are
independent of $x$, satisfied by the Sobolev-type Laguerre polynomials
$\set{\gL}$.

Further we show that in the case $M=0$ and $N>0$ the polynomials are
eigenfunctions of a linear differential operator, which is of order $2\a+8$
if $\a$ is a nonnegative integer and of infinite order otherwise.

Finally, we show that in the case $M>0$ and $N>0$ the polynomials are
eigenfunctions of a linear differential operator, which is of order $4\a+10$
if $\a$ is a nonnegative integer and of infinite order otherwise.
\end{abstract}

\vfill

\begin{tabular}{l}
Keywords : Differential equations, Sobolev-type Laguerre polynomials.\\
\\
1991 Mathematics Subject Classification : Primary 33C45 ; Secondary 34A35.
\end{tabular}

\section{Introduction}

Let ${\cal P}$ denote the space of all real polynomials. We consider the
polynomials $\set{\gL}$ which are orthogonal with respect to the
Sobolev-type inner product
$$<f,g>\;=\frac{1}{\G(\a+1)}\int\limits_0^{\infty}x^{\a}e^{-x}f(x)g(x)dx+
Mf(0)g(0)+Nf'(0)g'(0),\;f,g\in{\cal P},$$
where $\a>-1$, $M\ge 0$ and $N\ge 0$. These Sobolev-type Laguerre
polynomials (see \cite{Thesis} and \cite{SIAM}) are generalizations
of the generalized Laguerre polynomials $\set{\KL}$ found by
T.H.~Koornwinder in \cite{Koorn}. They can be written as
\be\la{def}\gL=A_0\L+A_1D\L+A_2D^2\L,\;\ndots,\ee
where $D=\frac{d}{dx}$ denotes the differentiation operator and the
coefficients $A_0$, $A_1$ and $A_2$ are given by
\be\la{AAA}\l\{\ba{l}\ds A_0=1+M{n+\a \ch n-1}+\frac{n(\a+2)-(\a+1)}{(\a+1)(\a+3)}
N{n+\a \ch n-2}+{}\\
\ds\hspace{5cm}{}+\frac{MN}{(\a+1)(\a+2)}{n+\a \ch n-1}{n+\a+1 \ch n-2}\\  \\
\ds A_1=M{n+\a \ch n}+\frac{n-1}{\a+1}N{n+\a \ch n-1}+
\frac{2MN}{(\a+1)^2}{n+\a \ch n}{n+\a+1 \ch n-2}\\  \\
\ds A_2=\frac{N}{\a+1}{n+\a \ch n-1}+\frac{MN}{(\a+1)^2}{n+\a \ch n}
{n+\a+1 \ch n-1}.\ea\r.\ee
For details concerning these Sobolev-type Laguerre polynomials and their
definition the reader is referred to \cite{Thesis} and \cite{SIAM}.

We consider differential equations of the form
\be\la{algdv}\sum_{i=1}^{\infty}d_i(x)y^{(i)}(x)=\lambda_ny(x),\ee
satisfied by these Sobolev-type Laguerre polynomials, where the coefficients
$\l\{d_i(x)\r\}_{i=1}^{\infty}$ are independent of $n$. For convenience we
write these differential equations in the form
$$\sum_{i=0}^{\infty}e_i(x)y^{(i)}(x)+xy''(x)+(\a+1-x)y'(x)+ny(x)=0,$$
where the coefficients $\l\{e_i(x)\r\}_{i=1}^{\infty}$ are independent of
$n$ and the coefficient $e_0(x):=e_0(n,\a)$ is independent of $x$.

The following lemma is well-known and easy to check. See for instance
\cite{KrallSheffer}~:

\lem{Let $\set{p_n(x)}$ be an arbitrary set of polynomials with
degree$[p_n(x)]=n$ for each $\ndots$ and let $\set{\lambda_n}$ be an
arbitrary sequence of constants with $\lambda_0=0$ and
$\l\{\lambda_n\r\}_{n=1}^{\infty}$ not all equal to zero. Then there exists
a unique sequence $\l\{d_i(x)\r\}_{i=1}^{\infty}$ of polynomials with
degree$[d_i(x)]\le i$ \footnote{The degree of the zero polynomial is
supposed to be $-1$} for all $i=1,2,3,\ldots$, such that (\ref{algdv})
has $\set{p_n(x)}$ as a polynomial set of solutions. Moreover, if
$d_i(x)=k_ix^i+$ lower order terms for $i=1,2,3,\ldots$, then
$$nk_1+n(n-1)k_2+\ldots+n!k_n=\lambda_n,\;n=1,2,3,\ldots.$$}

In \cite{DV} J.~Koekoek and R.~Koekoek found a differential equation of the
form (\ref{algdv}) for the generalized Laguerre polynomials $\set{\KL}$.
These polynomials form a special case of the above mentioned Sobolev-type
Laguerre polynomials, since $L_n^{\a,M,0}(x)=\KL$. For the symmetric form
of that differential equation the reader is referred to \cite{Sym}.

In \cite{Bav} H.~Bavinck found a new method to obtain the main result of
\cite{DV}. In this paper we will use this method to construct all
differential equations of the form
\bea\la{Sobdv}& &M\sum_{i=0}^{\infty}a_i(x)y^{(i)}(x)+
N\sum_{i=0}^{\infty}b_i(x)y^{(i)}(x)+{}\nn
& &\hspace{1cm}{}+MN\sum_{i=0}^{\infty}c_i(x)y^{(i)}(x)+
xy''(x)+(\a+1-x)y'(x)+ny(x)=0,\eea
where the coefficients $\l\{a_i(x)\r\}_{i=1}^{\infty}$,
$\l\{b_i(x)\r\}_{i=1}^{\infty}$ and $\l\{c_i(x)\r\}_{i=1}^{\infty}$ are
independent of $n$ and the coefficients $a_0(x)$, $b_0(x)$ and $c_0(x)$ are
independent of $x$, satisfied by the polynomials $\set{\gL}$ given by
(\ref{def}) and (\ref{AAA}).

In view of lemma~1 the coefficient $Ma_i(x)+Nb_i(x)+MNc_i(x)$ must be a
polynomial, independent of $n$, of degree at most $i$ for each
$i=1,2,3,\ldots$. Since $M\ge 0$ and $N\ge 0$ are arbitrary we conclude that
$\l\{a_i(x)\r\}_{i=1}^{\infty}$, $\l\{b_i(x)\r\}_{i=1}^{\infty}$ and
$\l\{c_i(x)\r\}_{i=1}^{\infty}$ must be polynomials, independent of $n$,
with degree$[a_i(x)]\le i$, degree$[b_i(x)]\le i$ and degree$[c_i(x)]\le i$
for every $i=1,2,3,\ldots$.

The main results of this paper given in the next section (see theorem~3)
give rise to the following corollary.

\cor{If $\a$ is a nonnegative integer, the polynomials $\set{\gL}$ satisfy a
differential equation of formal order $4\a+10$ which is of the form
\bea & &M\sum_{i=0}^{2\a+4}a_i(x)y^{(i)}(x)+
N\sum_{i=0}^{2\a+8}\b_i(x)y^{(i)}(x)+{}\nn
& &\hspace{1cm}{}+MN\sum_{i=0}^{4\a+10}\c_i(x)y^{(i)}(x)+
xy''(x)+(\a+1-x)y'(x)+ny(x)=0,\n\eea
where the coefficients $\l\{a_i(x)\r\}_{i=1}^{2\a+4}$,
$\l\{\b_i(x)\r\}_{i=1}^{2\a+8}$ and $\l\{\c_i(x)\r\}_{i=1}^{4\a+10}$ are
polynomials independent of $n$ which satisfy
$$\sum_{i=1}^{2\a+4}a_i(x)=\sum_{i=1}^{2\a+8}\b_i(x)=
\sum_{i=1}^{4\a+10}\c_i(x)=0.$$}
This corollary was stated as a conjecture in \cite{JCAM}. By formal order
we mean that for special cases ($M=0$ or $N=0$) the true order might be
lower. In fact, for $M>0$ and $N=0$ we find the order $2\a+4$ as in
\cite{DV} and for $M=0$ and $N>0$ we find a differential equation of order
$2\a+8$ when $\a$ is a nonnegative integer.

We emphasize that except for the examples given in \cite{JCAM} ($\a=0$) and
in \cite{Rap} ($\a=0$, $\a=1$ and $\a=2$) this is the first paper on
differential equations of the form (\ref{algdv}) satisfied by Sobolev-type
orthogonal polynomials.

For more results on Sobolev orthogonality and spectral differential
equations the reader is referred to \cite{Kwon1}. Some results concerning
the symmetrizability of the differential equations obtained in this paper
can be found in \cite{Kwon2}.

Finally, we refer to \cite{Charlier} and \cite{SobCharlier} where difference
operators are found for generalized Charlier polynomials and to
\cite{Meixner} where generalizations of Meixner polynomials are treated. In
these discrete cases the difference operators turn out to be of infinite
order for all values of the parameters.

The results obtained in \cite{DV} and in this paper generalize the fourth
order differential equation for the so-called Laguerre type polynomials
($\a=0$) found by H.L.~Krall in \cite{HLKrall1} (see also \cite{HLKrall2}).
These Laguerre type polynomials are described in more details in
\cite{AMKrall}.

\section{The main results}

We look for all differential equations of the form (\ref{Sobdv}) satisfied
by the Sobolev-type Laguerre polynomials $\set{\gL}$ defined by (\ref{def})
and (\ref{AAA}). We emphasize that we demand that the coefficients
$\l\{a_i(x)\r\}_{i=1}^{\infty}$, $\l\{b_i(x)\r\}_{i=1}^{\infty}$ and
$\l\{c_i(x)\r\}_{i=1}^{\infty}$ are independent of $n$ and that $a_0(x)$,
$b_0(x)$ and $c_0(x)$ do not depend on $x$. Sometimes we will use the
notation
$$a_0(x)=a_0(n,\a),\;b_0(x)=b_0(n,\a)\;\mbox{ and }\;c_0(x)=c_0(n,\a),
\;\ndots$$
in order to denote the dependence of $n$.

We will prove that the coefficients
$\l\{a_i(x)\r\}_{i=0}^{\infty}$ are given by
\be\la{anul}a_0(x)=a_0(n,\a)={n+\a+1 \ch n-1},\;\ndots\ee
and
\be\la{ai}a_i(x)=\frac{1}{i!}\sum_{j=1}^i(-1)^{i+j+1}
{\a+1 \ch j-1}{\a+2 \ch i-j}(\a+3)_{i-j}x^j,\;i=1,2,3,\ldots.\ee
This was already found in \cite{DV} and also in \cite{Bav}.

From (\ref{ai}) it is not very difficult to see that, for
$\a\ne 0,1,2,\ldots$
$$\mbox{degree}[a_i(x)]=i,\;i=1,2,3,\ldots$$
and for nonnegative integer values of $\a$
$$\l\{\ba{ll}\mbox{degree}[a_i(x)]=i, & i=1,2,3,\ldots,\a+2,\\ \\
\mbox{degree}[a_i(x)]=\a+2, & i=\a+3,\a+4,\a+5,\ldots,2\a+4,\\ \\
a_i(x)=0, & i>2\a+4.\ea\r.$$
Note that we have
$$a_{2\a+4}(x)=(-1)^{\a+1}\frac{x^{\a+2}}{(\a+2)!}\;\mbox{ for }\;
\a\in\{0,1,2,\ldots\}.$$
This implies that, for $M>0$ the differential equation given by
$$M\sum_{i=0}^{\infty}a_i(x)y^{(i)}(x)+xy''(x)+(\a+1-x)y'(x)+ny(x)=0,$$
satisfied by the polynomials $\set{\KL}$, where the coefficients
$\l\{a_i(x)\r\}_{i=0}^{\infty}$ are given by (\ref{anul}) and (\ref{ai}),
has order $2\a+4$ if $\a$ is a nonnegative integer and has infinite order
otherwise.

The coefficients $\l\{b_i(x)\r\}_{i=0}^{\infty}$ and
$\l\{c_i(x)\r\}_{i=0}^{\infty}$ are not unique. In fact we will show that
\be\la{bsplit}\l\{\ba{l}b_0(0,\a)=0\\ \\
b_0(n,\a)=b_0(1,\a)+\b_0(n,\a),\;n=1,2,3,\ldots\\ \\
b_i(x)=b_0(1,\a)b_i^*(x)+\b_i(x),\;i=1,2,3,\ldots\ea\r.\ee
and
\be\la{csplit}\l\{\ba{l}c_0(0,\a)=0\\ \\
c_0(n,\a)=b_0(1,\a)+\c_0(n,\a),\;n=1,2,3,\ldots\\ \\
c_i(x)=b_0(1,\a)c_i^*(x)+\c_i(x),\;i=1,2,3,\ldots,\ea\r.\ee
where $b_0(1,\a)$ is arbitrary. Here we have (cf. \cite{Rap} and
\cite{JCAM})
\be\la{bster}b_i^*(x)=\frac{1}{i!}
\sum_{j=0}^i(-1)^j{i \ch j}(\a+1)_{i-j}x^j,\;i=1,2,3,\ldots\ee
and
\be\la{cster}c_i^*(x)=\frac{(-1)^i}{i!}x^i,\;i=1,2,3,\ldots.\ee

Further we will show that the coefficients $\l\{\b_i(x)\r\}_{i=0}^{\infty}$
are given by
\be\la{bnul}\b_0(x)=\b_0(n,\a)=
\frac{n(\a+2)-\a}{(\a+1)(\a+4)}{n+\a+1 \ch n-2},\;n=1,2,3,\ldots,\ee
\be\la{bi0}\b_1(x)=0\;\mbox{ and }\;\b_i(x)=\b_i^{(1)}(x)+\b_i^{(2)}(x)+
\b_i^{(3)}(x)+\b_i^{(4)}(x),\;i=2,3,4,\ldots,\ee
where
\be\la{bi1}\b_i^{(1)}(x)=\frac{1}{(\a+1)(i-2)!}
\sum_{j=1}^i(-1)^{i+j+1}{\a+2 \ch j-1}{\a+3 \ch i-j}(\a+2)_{i-j}x^{j-1},
\;i=2,3,4,\ldots,\ee
\be\la{bi2}\b_i^{(2)}(x)=\frac{2(\a+2)}{(\a+1)(i-1)!}
\sum_{j=2}^i(-1)^{i+j+1}{\a+1 \ch j-2}{\a+2 \ch i-j}(\a+3)_{i-j}x^{j-2},
\;i=2,3,4,\ldots,\ee
\bea\la{bi3}& &\b_i^{(3)}(x)=\frac{1}{(\a+1)^2(\a+2)(\a+3)^2(i-2)!}
\sum_{j=0}^i(-1)^{i+j+1}{\a+3 \ch j}{\a+3 \ch i-j}(\a+1)_{i-j}\times{}\nn
& &{}\hspace{3cm}{}\times\l[(\a+1)(\a+3)+j(i-j)\r]x^j,\;i=2,3,4,\ldots\eea
and
\bea\la{bi4}& &\b_i^{(4)}(x)=\frac{i-2}{(\a+1)^2(\a+3)(\a+4)^2(i-2)!}
\sum_{j=0}^i(-1)^{i+j+1}{\a+4 \ch j}{\a+4 \ch i-j}(\a+1)_{i-j}\times{}\nn
& &{}\hspace{3cm}{}\times\l[(\a+1)(\a+4)+j(i-j)\r]x^j,\;i=2,3,4,\ldots.\eea
From (\ref{bi0}), (\ref{bi1}), (\ref{bi2}), (\ref{bi3}) and (\ref{bi4}) we
obtain that for $\a\ne 0,1,2,\ldots$
$$\mbox{degree}[\b_i(x)]=i,\;i=2,3,4,\ldots$$
and for nonnegative integer values of $\a$
$$\l\{\ba{ll}\mbox{degree}[\b_i(x)]=i, & i=2,3,4,\ldots,\a+4,\\ \\
\mbox{degree}[\b_i(x)]=\a+4, & i=\a+5,\a+6,\a+7,\ldots,2\a+8,\\ \\
\b_i(x)=0, & i>2\a+8.\ea\r.$$
We remark that
$$\b_{2\a+8}(x)=\b_{2\a+8}^{(4)}(x)=
(-1)^{\a+1}\frac{\a+2}{\a+1}\frac{x^{\a+4}}{(\a+4)!}\;\mbox{ for }\;
\a\in\{0,1,2,\ldots\}.$$
This implies that for $N>0$ the differential equation
$$N\sum_{i=0}^{\infty}\b_i(x)y^{(i)}(x)+xy''(x)+(\a+1-x)y'(x)+ny(x)=0,$$
satisfied by the polynomials $\set{L_n^{\a,0,N}(x)}$, where the coefficients
$\l\{\b_i(x)\r\}_{i=0}^{\infty}$ are given by $\b_0(0,\a)=0$ and
(\ref{bnul}), (\ref{bi0}), (\ref{bi1}), (\ref{bi2}), (\ref{bi3}) and
(\ref{bi4}), has order $2\a+8$ if $\a$ is a nonnegative integer and has
infinite order otherwise.

Finally we will show that the coefficients $\l\{\c_i(x)\r\}_{i=0}^{\infty}$
are given by
\be\la{cnul}\c_0(x)=\c_0(n,\a)=\frac{1}{\a+1}{n+\a+1 \ch n-2}
\hyp{3}{2}{-n+2,-\a-2,\a+3}{2,\a+4}{1},\;n=1,2,3,\ldots\ee
and
\be\la{ci0}\c_i(x)=\c_i^{(1)}(x)+\c_i^{(2)}(x)+\c_i^{(3)}(x)+
\c_i^{(4)}(x)+\c_i^{(5)}(x),\;i=1,2,3,\ldots,\ee
where for $i=1,2,3,\ldots$
\be\la{ci1}\c_i^{(1)}(x)=\frac{(-1)^{i+1}}{(\a+1)(i-1)!}\sum_{k=1}^ia_k(x)
\sum_{j=k}^i(-1)^j{\a+3 \ch j-k}{\a+3 \ch i-j}(\a+4)_{i-j+k-1}x^{j-k},\ee
\be\la{ci2}\c_i^{(2)}(x)=\frac{(-1)^i}{(i+1)!}\sum_{k=1}^i\b_k(x)
\sum_{j=k}^i(-1)^j{\a+1 \ch j-k}{\a+1 \ch i-j}(\a+3)_{i-j+k}x^{j-k},\ee
\bea\la{ci3}& &\c_i^{(3)}(x)=\frac{2}{(\a+1)(\a+3)}\sum_{j=0}^i\sum_{n=0}^j
(-1)^j\frac{(-j)_n(\a+3)_{i-j+n}(\a+3)_n}{j!(i-j)!n!}\times{}\nn
& &{}\hspace{5cm}{}\times\sum_{k=0}^{i-j}
\frac{(-i+j)_k(n+\a+3)_k(n+\a+4)_k}{\G(n+k)\G(n+k+2)k!}x^j,\eea
\bea\la{ci4}& &\c_i^{(4)}(x)=-\frac{2}{(\a+1)(\a+4)}\sum_{j=0}^i\sum_{n=0}^j
(-1)^j\frac{(-j)_n(\a+3)_{i-j+n}(\a+4)_n}{j!(i-j)!n!}\times{}\nn
& &{}\hspace{5cm}{}\times\sum_{k=0}^{i-j}
\frac{(-i+j)_k(n+\a+3)_k(n+\a+4)_k}{\G(n+k)\G(n+k+2)k!}x^j\eea
and
\bea\la{ci5}& &\c_i^{(5)}(x)=\frac{(-1)^{i+1}}{(\a+1)(\a+3)}
\sum_{k=0}^{i-1}(-1)^k\frac{(-\a-2)_k(\a+3)_k}{(2)_k(\a+4)_kk!}\times{}\nn
& &{}\hspace{3cm}{}\times\frac{1}{\G(i-k)}
\sum_{j=0}^i(-1)^j{k+\a+3 \ch j}{k+\a+3 \ch i-j}(\a+3)_{i-j}x^j.\eea
From (\ref{ci0}), (\ref{ci1}), (\ref{ci2}), (\ref{ci3}), (\ref{ci4}) and
(\ref{ci5}) we obtain that $\c_1(x)=0$.

Moreover we will show that, for nonnegative integer values of $\a$
$$\c_i(x)=0\;\mbox{ for }\;i>4\a+10$$
and that
$$\c_{4\a+10}(x)=\c_{4\a+10}^{(5)}(x)=
\frac{x^{2\a+5}}{(\a+1)(2\a+5)(\a+2)!(\a+3)!}\;\mbox{ for }\;
\a\in\{0,1,2,\ldots\}.$$

Note that for $i\ge 3$ the coefficient of $y^{(i)}(x)$ in the differential
equation (\ref{Sobdv}) is equal to $Ma_i(x)+Nb_i(x)+MNc_i(x)$. Hence for
$M^2+N^2>0$ finite order can only occur when $Ma_i(x)$, $Nb_i(x)$ and
$MNc_i(x)$ all vanish for $i\ge K$ for some $K\in\{3,4,5,\ldots\}$.

Now we have reached our main theorem~:

\st{For $\a>-1$ and $M^2+N^2>0$ the only differential equations of the form
(\ref{Sobdv}) satisfied by the Sobolev-type Laguerre polynomials
$\set{\gL}$ defined by (\ref{def}) and (\ref{AAA}), are those where the
coefficients $\l\{a_i(x)\r\}_{i=0}^{\infty}$,
$\l\{b_i(x)\r\}_{i=0}^{\infty}$ and $\l\{c_i(x)\r\}_{i=0}^{\infty}$ are
given by (\ref{anul}) up to and including (\ref{ci5}) with $b_0(1,\a)$
arbitrary.

Only if $Nb_0(1,\a)=0$ and $\a\in\{0,1,2,\ldots\}$ the order of this
differential equation is finite and equal to
$$\l\{\ba{lcl}2\a+4 & \mbox{if} & M>0\;\mbox{ and }\;N=0\\
2\a+8 & \mbox{if} & M=0\;\mbox{ and }\;N>0\\
4\a+10 & \mbox{if} & M>0\;\mbox{ and }\;N>0.\ea\r.$$
Otherwise the differential equation is of infinite order.}

Finally we will show that the coefficients $\l\{a_i(x)\r\}_{i=1}^{\infty}$,
$\l\{\b_i(x)\r\}_{i=1}^{\infty}$ and $\l\{\c_i(x)\r\}_{i=1}^{\infty}$
satisfy
$$\sum_{i=1}^{\infty}a_i(x)=-\frac{\sin\pi\a}{\pi}\frac{x}{(\a+2)(\a+3)}
\hyp{1}{1}{1}{\a+4}{-x},\;\a>-1,$$
$$\sum_{i=1}^{\infty}\b_i(x)=\frac{\sin\pi\a}{\pi}
\frac{\a+2}{(\a+1)(\a+3)(\a+4)}\l[1-x\,\hyp{2}{2}{1,\a+3}{\a+2,\a+5}{-x}\r],
\;\a>-1$$
and
$$\sum_{i=1}^{\infty}\c_i(x)=0\;\mbox{ for }\;\a\in\{0,1,2,\ldots\}.$$
This implies that for nonnegative integer values of $\a$ we have
$$\sum_{i=1}^{2\a+4}a_i(x)=\sum_{i=1}^{2\a+8}\b_i(x)=
\sum_{i=1}^{4\a+10}\c_i(x)=0.$$

\section{Some classical formulas}

We start with the following lemma on partial sums of Gauss' hypergeometric
series which will be used in this paper~:

\lem{Let $n\in\{0,1,2,\ldots\}$. Then the $n$th partial sum of the Gauss'
hypergeometric series at the point 1
$$\sum_{k=0}^{\infty}\frac{(a)_k(b)_k}{(c)_kk!},$$
equals
$$\sum_{k=0}^n\frac{(a)_k(b)_k}{(c)_kk!}=
\frac{(a+1)_n}{n!}\hyp{3}{2}{-n,a,c-b}{a+1,c}{1}$$
for all values of the parameters $a$, $b$ and $c$ for which
$(a+1)_n(c)_n\ne 0$.}

{\bf Proof.} The proof is based on the well-known Vandermonde summation
formula
\be\la{Van}\hyp{2}{1}{-n,b}{c}{1}=\frac{(c-b)_n}{(c)_n},\;(c)_n\ne 0,
\;\ndots.\ee
Suppose that $n\in\{0,1,2,\ldots\}$. Then we find for $(a+1)_n(c)_n\ne 0$
by changing the order of summation and applying (\ref{Van}) twice
\bea & &\hyp{3}{2}{-n,a,c-b}{a+1,c}{1}=
\sum_{m=0}^n\frac{(-n)_m(a)_m}{(a+1)_mm!}\hyp{2}{1}{-m,b}{c}{1}\nn
&=&\sum_{m=0}^n\sum_{k=0}^m
\frac{(-n)_m(a)_m}{(a+1)_mm!}\frac{(-m)_k(b)_k}{(c)_kk!}=
\sum_{k=0}^n\sum_{m=k}^n
\frac{(-n)_m(a)_m(-m)_k(b)_k}{(a+1)_m(c)_km!k!}\nn
&=&\sum_{k=0}^n\sum_{m=0}^{n-k}
\frac{(-n)_{k+m}(a)_{k+m}(-k-m)_k(b)_k}{(a+1)_{k+m}(c)_k(k+m)!k!}\nn
&=&\sum_{k=0}^n\frac{(-n)_k(a)_k(b)_k}{(a+1)_k(c)_kk!}
(-1)^k\hyp{2}{1}{-n+k,a+k}{a+k+1}{1}\nn
&=&\sum_{k=0}^n\frac{(-n)_k(a)_k(b)_k}{(a+1)_k(c)_kk!}
(-1)^k\frac{(1)_{n-k}}{(a+k+1)_{n-k}}=
\frac{n!}{(\a+1)_n}\sum_{k=0}^n\frac{(a)_k(b)_k}{(c)_kk!}.\n\eea
This proves lemma~4.

The special case $b=c$ of lemma~4 leads to the $n$th partial sum of a
${}_1F_0$ hypergeometric series at the point $1$~:
\be\la{binform}\sum_{k=0}^n{a+k \ch k}=\sum_{k=0}^n\frac{(a+1)_k}{k!}=
\frac{(a+2)_n}{n!}={n+a+1 \ch n},\;\ndots.\ee

Further we list some definitions and properties of the classical Laguerre
polynomials $\set{\L}$ which we need in this paper.

For $\a$ real and $\a>-1$ the classical Laguerre polynomials are orthogonal
on the interval $[0,\infty)$ with respect to the weight function
$x^{\a}e^{-x}$. They are usually defined by
\be\la{defhyp}\L={n+\a \ch n}\hyp{1}{1}{-n}{\a+1}{x},\;\ndots.\ee

Note that $\L$ is also a polynomial in the parameter $\a$. This implies that
the classical Laguerre polynomials can be defined for all $\a$ by
\bea\la{def1}\L
&=&\frac{1}{n!}\sum_{k=0}^n(-n)_k(\a+k+1)_{n-k}\frac{x^k}{k!}\\
&=&\la{def2}(-1)^n\sum_{k=0}^n\frac{(-n-\a)_{n-k}}{(n-k)!}\frac{x^k}{k!}\\
&=&\la{def3}\sum_{k=0}^n(-1)^k\frac{(\a+k+1)_{n-k}}{(n-k)!}\frac{x^k}{k!},
\;\ndots.\eea

Now we take $\a$ arbitrary. We have
\be\la{nul}L_n^{(\a)}(0)={n+\a \ch n}=\frac{(\a+1)_n}{n!},\;\ndots\ee
and
\be\la{diff}D^k\L=(-1)^kL_{n-k}^{(\a+k)}(x),\;k\le n,\;k,\ndots.\ee

The classical Laguerre polynomials satisfy the second order linear
differential equation
\be\la{dvLag}xy''(x)+(\a+1-x)y'(x)+ny(x)=0.\ee

From the generating function given by
$$\sum_{n=0}^{\infty}\L t^n=(1-t)^{-\a-1}\exp\l(\frac{xt}{t-1}\r),$$
it easily follows that for arbitrary $p$ we have
\be\la{rel}L_n^{(\a-p)}(x)=\sum_{k=0}^n(-1)^k{p \ch k}L_{n-k}^{(\a)}(x),
\;\ndots.\ee
Further we have for $i=0,1,2,\ldots$
\be\la{rel1}D^i\L=D^{i+1}\L-D^{i+1}\Ln,\;\ndots,\ee
or equivalently
\be\la{rel2}D^{i+1}\Ln=D^{i+1}\L-D^i\L,\;\ndots.\ee

We also obtain from the generating function that
$$\sum_{k=0}^{\infty}L_k^{(-\a-i-1)}(-x)t^k
\sum_{m=0}^{\infty}L_m^{(\a+j)}(x)t^m=(1-t)^{i-j-1}.$$
This implies that
\be\la{inv}\sum_{k=j}^iL_{i-k}^{(-\a-i-1)}(-x)L_{k-j}^{(\a+j)}(x)=
\delta_{ij},\;j\le i,\;i,j=0,1,2,\ldots.\ee
This inversion formula is an important tool to obtain the results of this
paper. We will use it in the following way~:

\lem{Suppose that for $k\in\{0,1,2,\ldots\}$ we have the system of equations
$$\sum_{i=1}^{\infty}A_i(x)D^{i+k}\L=F_n(x),\;n=k+1,k+2,k+3,\ldots,$$
where $\l\{A_i(x)\r\}_{i=1}^{\infty}$ are independent of $n$.
Then this system has a unique solution given by
$$A_i(x)=(-1)^{i+k}\sum_{j=1}^iL_{i-j}^{(-\a-i-k-1)}(-x)F_{j+k}(x),
\;i=1,2,3,\ldots.$$}

Finally we define the following polynomials involving the classical Laguerre
polynomials
\be\la{F}F_i(a,b;c;x):=\sum_{k=0}^i\frac{(a)_k}{\G(c+k)}
L_{i-k}^{(-a-i)}(-x)L_k^{(b)}(x),\;i=0,1,2,\ldots,\ee
\be\la{G}G_i(a,b;c;x):=\sum_{k=0}^i\frac{(a)_k}{\G(c+k)}
L_{i-k}^{(-b-i-1)}(-x)L_k^{(b)}(x),\;i=0,1,2,\ldots\ee
and
\be\la{H}H_i(a,b,c;d,e;x):=\sum_{k=0}^i\frac{(a)_k(b)_k}{\G(d+k)\G(e+k)}
L_{i-k}^{(-a-i)}(-x)L_k^{(c)}(x),\;i=0,1,2,\ldots.\ee
Note that $F_i$ and $G_i$ are polynomials in $a$, $b$ and $x$ and that $H_i$
is a polynomial in $a$, $b$, $c$ and $x$ for each $i=0,1,2,\ldots$.

Now we will prove the following lemma~:

\lem{
\be\la{lemF}F_i(a,b;c;x)=\frac{1}{\G(c+i)}
\sum_{j=0}^i(-1)^j{a-c \ch j}{b-c+1 \ch i-j}(a)_{i-j}x^j,\;i=0,1,2,\ldots\ee
and
\bea\la{lemG}& &G_i(a+2,a;c;x)\nn
&=&\frac{1}{(a+1)(a-c+2)\G(c+i)}
\sum_{j=0}^i(-1)^j{a-c+2 \ch j}{a-c+2 \ch i-j}(a+1)_{i-j}\times{}\nn
& &{}\hspace{3cm}{}\times\l[(a+1)(a-c+2)+j(i-j)\r]x^j,
\;i=0,1,2,\ldots.\eea}

{\bf Proof.} First we define (see for instance \cite{Meijer})
$$\hyphyp{2}{1}{a,b}{c}{z}:=
\sum_{k=0}^{\infty}\frac{(a)_k(b)_k}{\G(c+k)k!}z^k,$$
for all complex $a$, $b$, $c$ and $z$ for which the series in the right-hand
side converges. If $c\ne 0,-1,-2,\ldots$ we may write
$$\hyphyp{2}{1}{a,b}{c}{z}=\frac{1}{\G(c)}\hyp{2}{1}{a,b}{c}{z}.$$
Now we have (cf. (\ref{Van}))
\be\la{som}\hyphyp{2}{1}{-n,b}{c}{1}=\frac{(c-b)_n}{\G(c+n)},\;\ndots,\ee
for all complex values of $b$ and $c$. This is the well-known Vandermonde
summation formula.

If we apply definition (\ref{def2}) to $L_{i-k}^{(-a-i)}(-x)$ and definition
(\ref{def3}) to $L_k^{(b)}(x)$ in (\ref{F}), change the order of the
summations and apply (\ref{som}) twice we obtain for each $i=0,1,2,\ldots$
\bea F_i(a,b;c;x)&=&\sum_{k=0}^i\sum_{m=0}^{i-k}\sum_{n=0}^k(-1)^{i+n}
\frac{(a)_k(-i)_{k+m}(a+k)_{i-k-m}(n+b+1)_{k-n}}{\G(c+k)i!m!(k-n)!n!}x^{m+n}\nn
&=&\sum_{m=0}^i\sum_{n=0}^{i-m}\sum_{k=0}^{i-m-n}(-1)^{i+n}
\frac{(-i)_{m+n+k}(a)_{i-m}(n+b+1)_k}{\G(c+n+k)i!m!n!k!}x^{m+n}\nn
&=&\sum_{m=0}^i\sum_{n=0}^{i-m}(-1)^{i+n}\frac{(-i)_{m+n}(a)_{i-m}}{i!m!n!}
\hyphyp{2}{1}{-i+m+n,n+b+1}{c+n}{1}x^{m+n}\nn
&=&\sum_{j=0}^i\sum_{n=0}^j(-1)^{i+n}\frac{(-i)_j(a)_{i-j+n}(c-b-1)_{i-j}}
{i!(j-n)!n!\G(c+i-j+n)}x^j\nn
&=&\frac{(-1)^i}{i!}\sum_{j=0}^i\frac{(-i)_j(a)_{i-j}(c-b-1)_{i-j}}{j!}
\hyphyp{2}{1}{-j,a+i-j}{c+i-j}{1}x^j\nn
&=&\frac{(-1)^i}{i!\G(c+i)}\sum_{j=0}^i(-i)_j(c-a)_j(a)_{i-j}(c-b-1)_{i-j}
\frac{x^j}{j!}\nn
&=&\frac{1}{\G(c+i)}
\sum_{j=0}^i(-1)^j{a-c \ch j}{b-c+1 \ch i-j}(a)_{i-j}x^j,\n\eea
which proves (\ref{lemF}).

In a similar way we obtain for each $i=0,1,2,\ldots$
\bea G_i(a+2,a;c;x)&=&\sum_{k=0}^i\sum_{m=0}^{i-k}\sum_{n=0}^k(-1)^{i+n}
\frac{(-i)_{k+m}(a+2)_k(n+a+1)_{i-m-n}}{\G(c+k)i!m!(k-n)!n!}x^{m+n}\nn
&=&\sum_{m=0}^i\sum_{n=0}^{i-m}\sum_{k=0}^{i-m-n}(-1)^{i+n}
\frac{(-i)_{m+n+k}(a+2)_{n+k}(n+a+1)_{i-m-n}}{\G(c+n+k)i!m!n!k!}x^{m+n}\nn
&=&\sum_{m=0}^i\sum_{n=0}^{i-m}(-1)^{i+n}
\frac{(-i)_{m+n}(a+2)_n(n+a+1)_{i-m-n}}{i!m!n!}\times{}\nn
& &\hspace{4cm}{}\times\hyphyp{2}{1}{-i+m+n,n+a+2}{n+c}{1}x^{m+n}\nn
&=&\sum_{j=0}^i\sum_{n=0}^j(-1)^{i+n}
\frac{(-i)_j(a+2)_n(n+a+1)_{i-j}(c-a-2)_{i-j}}{i!(j-n)!n!\G(c+i-j+n)}x^j\nn
&=&\frac{(-1)^i}{i!}\sum_{j=0}^i(-i)_j(c-a-2)_{i-j}\frac{x^j}{j!}
\sum_{n=0}^j\frac{(-j)_n(a+2)_n(n+a+1)_{i-j}}{n!\G(c+i-j+n)}.\n\eea
Now we use (\ref{som}) to find for $j=1,2,3,\ldots$
\bea & &\sum_{n=0}^j\frac{(-j)_n(a+2)_n(n+a+1)_{i-j}}{n!\G(c+i-j+n)}
=\sum_{n=0}^j\frac{(-j)_n(a+1+n)(a+1)_{i-j+n}}{(a+1)n!\G(c+i-j+n)}\nn
&=&(a+1)_{i-j}\sum_{n=0}^j\frac{(-j)_n(a+i-j+1)_n}{n!\G(c+i-j+n)}+
(a+2)_{i-j}\sum_{n=1}^j\frac{(-j)_n(a+i-j+2)_{n-1}}{(n-1)!\G(c+i-j+n)}\nn
&=&(a+1)_{i-j}\hyphyp{2}{1}{-j,a+i-j+1}{c+i-j}{1}-
j(a+2)_{i-j}\hyphyp{2}{1}{-j+1,a+i-j+2}{c+i-j+1}{1}\nn
&=&\frac{(a+1)_{i-j}(c-a-2)_j}{(a+1)(a-c+2)\G(c+i)}
\l[(a+1)(a-c+2-j)+j(a+1+i-j)\r].\n\eea
Note that this result is also true for $j=0$. Hence, we obtain for
$i=0,1,2,\ldots$
\bea & &G_i(a+2,a;c;x)\nn
&=&\frac{(-1)^i}{(a+1)(a-c+2)i!\G(c+i)}
\sum_{j=0}^i(-i)_j(c-a-2)_j(c-a-2)_{i-j}(a+1)_{i-j}\times{}\nn
& &{}\hspace{3cm}{}\times\l[(a+1)(a-c+2-j)+j(a+1+i-j)\r]\frac{x^j}{j!}\nn
&=&\frac{1}{(a+1)(a-c+2)\G(c+i)}
\sum_{j=0}^i(-1)^j{a-c+2 \ch j}{a-c+2 \ch i-j}(a+1)_{i-j}\times{}\nn
& &\hspace{5cm}{}\times\l[(a+1)(a-c+2)+j(i-j)\r]x^j,\n\eea
which proves (\ref{lemG}).

If we also define (cf. \cite{Meijer})
\be\la{defphi}\hyphyp{3}{2}{a,b,c}{d,e}{z}:=
\sum_{k=0}^{\infty}\frac{(a)_k(b)_k(c)_k}{\G(d+k)\G(e+k)k!}z^k,\ee
for all complex $a$, $b$, $c$, $d$, $e$ and $z$ for which the series in the
right-hand side converges, then we also have~:

\lem{
\bea\la{lemH}H_i(a,b,c;d,e;x)&=&\sum_{j=0}^i\sum_{n=0}^j(-1)^{i+j}
\frac{(-j)_n(a)_{i-j+n}(b)_n}{j!(i-j)!n!}\times{}\nn
& &{}\hspace{1cm}{}\times\hyphyp{3}{2}{-i+j,n+b,n+c+1}{n+d,n+e}{1}x^j,
\;i=0,1,2,\ldots.\eea}

{\bf Proof.} As before we apply definition (\ref{def2}) to
$L_{i-k}^{(-a-i)}(-x)$ and definition (\ref{def3}) to $L_k^{(c)}(x)$ in
(\ref{H}) and change the order of the summations to obtain for each
$i=0,1,2,\ldots$
\bea & &H_i(a,b,c;d,e;x)\nn
&=&\sum_{k=0}^i\sum_{m=0}^{i-k}\sum_{n=0}^k(-1)^{i+n}
\frac{(-i)_{k+m}(a)_{i-m}(b)_k(n+c+1)_{k-n}}
{\G(d+k)\G(e+k)i!m!(k-n)!n!}x^{m+n}\nn
&=&\sum_{m=0}^i\sum_{n=0}^{i-m}\sum_{k=0}^{i-m-n}(-1)^{i+n}
\frac{(-i)_{m+n+k}(a)_{i-m}(b)_{n+k}(n+c+1)_k}
{\G(d+n+k)\G(e+n+k)i!m!n!k!}x^{m+n}\nn
&=&\sum_{m=0}^i\sum_{n=0}^{i-m}(-1)^{i+n}\frac{(-i)_{m+n}(a)_{i-m}(b)_n}
{i!m!n!}\hyphyp{3}{2}{-i+m+n,n+b,n+c+1}{n+d,n+e}{1}x^{m+n}\nn
&=&\sum_{j=0}^i\sum_{n=0}^j(-1)^{i+n}\frac{(-i)_j(a)_{i-j+n}(b)_n}
{i!(j-n)!n!}\hyphyp{3}{2}{-i+j,n+b,n+c+1}{n+d,n+e}{1}x^j\nn
&=&\sum_{j=0}^i\sum_{n=0}^j(-1)^{i+j}\frac{(-j)_n(a)_{i-j+n}(b)_n}
{j!(i-j)!n!}\hyphyp{3}{2}{-i+j,n+b,n+c+1}{n+d,n+e}{1}x^j,\n\eea
which proves (\ref{lemH}).

\section{The computation of the coefficients}

We take $\a$ real and $\a>-1$.

If we substitute $y(x)=\gL$ in the differential equation (\ref{Sobdv}) and
we use (\ref{def}) and (\ref{dvLag}) we find~:

\bea & &MA_0\sum_{i=0}^{\infty}a_i(x)D^i\L+
MA_1\sum_{i=0}^{\infty}a_i(x)D^{i+1}\L+
MA_2\sum_{i=0}^{\infty}a_i(x)D^{i+2}\L+{}\nn
& &\hspace{5mm}{}+NA_0\sum_{i=0}^{\infty}b_i(x)D^i\L+
NA_1\sum_{i=0}^{\infty}b_i(x)D^{i+1}\L+
NA_2\sum_{i=0}^{\infty}b_i(x)D^{i+2}\L+{}\nn
& &\hspace{5mm}{}+MNA_0\sum_{i=0}^{\infty}c_i(x)D^i\L+
MNA_1\sum_{i=0}^{\infty}c_i(x)D^{i+1}\L+{}\nn
& &\hspace{1cm}{}+MNA_2\sum_{i=0}^{\infty}c_i(x)D^{i+2}\L+{}\nn
& &\hspace{1cm}{}+A_1\l[xD^3\L+(\a+1-x)D^2\L+nD\L\r]+{}\nn
& &\hspace{1cm}{}+A_2\l[xD^4\L+(\a+1-x)D^3\L+nD^2\L\r]=0,\;\ndots.\n\eea
If we interpret the left-hand side as a polynomial in $M$ and $N$ and we use
$$xD^3\L+(\a+1-x)D^2\L+nD\L=D\L-D^2\L=-D^2\Ln$$
and
$$xD^4\L+(\a+1-x)D^3\L+nD^2\L=2D^2\L-2D^3\L=-2D^3\Ln,$$
which follow from the differential equation (\ref{dvLag}) for the classical
Laguerre polynomials and (\ref{rel2}), and the definition (\ref{AAA}) of the
coefficients $A_0$, $A_1$ and $A_2$ we find eight systems of equations for
the coefficients $\l\{a_i(x)\r\}_{i=0}^{\infty}$,
$\l\{b_i(x)\r\}_{i=0}^{\infty}$ and $\l\{c_i(x)\r\}_{i=0}^{\infty}$. These
can be written as follows~:

$$\ba{lcllcl}M & : & S_1=0 \hspace{5cm} & N & : & S_5=0\\
M^2 & : & S_2=0 & N^2 & : & S_6=0\\
MN & : & S_3=0 & MN^2 & : & S_7=0\\
M^2N & : & S_4=0 & M^2N^2 & : & S_8=0\ea$$
for $\ndots$, where
\bea & &S_1:=\sum_{i=0}^{\infty}a_i(x)D^i\L-{n+\a \ch n}D^2\Ln,\nn
& &S_2:={n+\a \ch n-1}\sum_{i=0}^{\infty}a_i(x)D^i\L+
{n+\a \ch n}\sum_{i=0}^{\infty}a_i(x)D^{i+1}\L,\nn
& &S_3:=\frac{n(\a+2)-(\a+1)}{(\a+1)(\a+3)}{n+\a \ch n-2}\sum_{i=0}^{\infty}
a_i(x)D^i\L+{}\nn
& &\hspace{15mm}{}+\frac{n-1}{\a+1}{n+\a \ch n-1}
\sum_{i=0}^{\infty}a_i(x)D^{i+1}\L+{}\nn
& &\hspace{15mm}{}+\frac{1}{\a+1}{n+\a \ch n-1}\sum_{i=0}^{\infty}a_i(x)D^{i+2}\L
+{n+\a \ch n-1}\sum_{i=0}^{\infty}b_i(x)D^i\L+{}\nn
& &\hspace{15mm}{}+{n+\a \ch n}\sum_{i=0}^{\infty}b_i(x)D^{i+1}\L
+\sum_{i=0}^{\infty}c_i(x)D^i\L+{}\nn
& &\hspace{15mm}{}-\frac{2}{(\a+1)^2}{n+\a \ch n}
\l[{n+\a+1 \ch n-2}D^2\Ln+{n+\a+1 \ch n-1}D^3\Ln\r],\nn
& &S_4:=\frac{1}{(\a+1)(\a+2)}{n+\a \ch n-1}{n+\a+1 \ch n-2}\sum_{i=0}^{\infty}
a_i(x)D^i\L+{}\nn
& &\hspace{15mm}{}+\frac{2}{(\a+1)^2}{n+\a \ch n}{n+\a+1 \ch n-2}
\sum_{i=0}^{\infty}a_i(x)D^{i+1}\L+{}\nn
& &\hspace{15mm}{}+\frac{1}{(\a+1)^2}{n+\a \ch n}{n+\a+1 \ch n-1}
\sum_{i=0}^{\infty}a_i(x)D^{i+2}\L+{}\nn
& &\hspace{15mm}{}+{n+\a \ch n-1}\sum_{i=0}^{\infty}c_i(x)D^i\L
+{n+\a \ch n}\sum_{i=0}^{\infty}c_i(x)D^{i+1}\L,\nn
& &S_5:=\sum_{i=0}^{\infty}b_i(x)D^i\L-\frac{1}{\a+1}{n+\a \ch n-1}
\l[(n-1)D^2\Ln+2D^3\Ln\r],\nn
& &S_6:=\frac{n(\a+2)-(\a+1)}{(\a+1)(\a+3)}{n+\a \ch n-2}\sum_{i=0}^{\infty}
b_i(x)D^i\L+{}\nn
& &\hspace{15mm}{}+\frac{n-1}{\a+1}{n+\a \ch n-1}\sum_{i=0}^{\infty}
b_i(x)D^{i+1}\L+{}\nn
& &\hspace{15mm}{}+\frac{1}{\a+1}{n+\a \ch n-1}\sum_{i=0}^{\infty}
b_i(x)D^{i+2}\L,\nn
& &S_7:=\frac{1}{(\a+1)(\a+2)}{n+\a \ch n-1}{n+\a+1 \ch n-2}\sum_{i=0}^{\infty}
b_i(x)D^i\L+{}\nn
& &\hspace{15mm}{}+\frac{2}{(\a+1)^2}{n+\a \ch n}{n+\a+1 \ch n-2}
\sum_{i=0}^{\infty}b_i(x)D^{i+1}\L+{}\nn
& &\hspace{15mm}{}+\frac{1}{(\a+1)^2}{n+\a \ch n}{n+\a+1 \ch n-1}
\sum_{i=0}^{\infty}b_i(x)D^{i+2}\L+{}\nn
& &\hspace{15mm}{}+\frac{n(\a+2)-(\a+1)}{(\a+1)(\a+3)}{n+\a \ch n-2}
\sum_{i=0}^{\infty}c_i(x)D^i\L+{}\nn
& &\hspace{15mm}{}+\frac{n-1}{\a+1}{n+\a \ch n-1}\sum_{i=0}^{\infty}
c_i(x)D^{i+1}\L+{}\nn
& &\hspace{15mm}{}+\frac{1}{\a+1}{n+\a \ch n-1}
\sum_{i=0}^{\infty}c_i(x)D^{i+2}\L\n\eea
and
\bea & &S_8:=\frac{1}{(\a+1)(\a+2)}{n+\a \ch n-1}{n+\a+1 \ch n-2}
\sum_{i=0}^{\infty}c_i(x)D^i\L+{}\nn
& &\hspace{15mm}{}+\frac{2}{(\a+1)^2}{n+\a \ch n}{n+\a+1 \ch n-2}
\sum_{i=0}^{\infty}c_i(x)D^{i+1}\L+{}\nn
& &\hspace{15mm}{}+\frac{1}{(\a+1)^2}{n+\a \ch n}{n+\a+1 \ch n-1}
\sum_{i=0}^{\infty}c_i(x)D^{i+2}\L.\n\eea
The systems of equations $S_1=0$ and $S_2=0$ for $\ndots$ lead to the
solution for the coefficients $\l\{a_i(x)\r\}_{i=0}^{\infty}$ which was
already found in \cite{DV}. The systems of equations $S_5=0$ and $S_6=0$
for $\ndots$ will lead to the solution for the coefficients
$\l\{b_i(x)\r\}_{i=0}^{\infty}$ while the other four will eventually lead
to the solution for the coefficients $\l\{c_i(x)\r\}_{i=0}^{\infty}$.

\subsection{The computation of the coefficients
$\l\{a_i(x)\r\}_{i=0}^{\infty}$}

Since $\a>-1$ we have ${n+\a \ch n}\ne 0$ for all $\ndots$. Then the
systems of equations $S_1=0$ and $S_2=0$ for $\ndots$ are equivalent to
\be\la{a1}\sum_{i=0}^{\infty}a_i(x)D^i\L={n+\a \ch n}D^2\Ln,\;\ndots\ee
and
\be\la{a2}\sum_{i=0}^{\infty}a_i(x)D^{i+1}\L=-{n+\a \ch n-1}D^2\Ln,\;
\ndots.\ee
By considering equation (\ref{a1}) for $n=0$ and equation (\ref{a2}) for
$n=0$ and $n=1$ we conclude that $a_0(0,\a)=0$ and $a_0(1,\a)=1$. Since
$a_i(x)$ must be a polynomial in $x$ of degree at most $i$ for each
$i=1,2,3,\ldots$, we may write
$$a_i(x)=k_ix^i+\mbox{lower order terms},\;i=1,2,3,\ldots.$$
By comparing the coefficients of highest degree in (\ref{a1}) and
(\ref{a2}) we find~:
$$\frac{a_0(n,\a)}{n!}+\sum_{i=1}^n\frac{k_i}{(n-i)!}=0,\;n=1,2,3,\ldots$$
and
$$\frac{a_0(n,\a)}{(n-1)!}+\sum_{i=1}^{n-1}\frac{k_i}{(n-i-1)!}=
{n+\a \ch n-1}\frac{1}{(n-1)!},\;n=2,3,4,\ldots.$$
Hence, since $k_i$ is independent of $n$, we obtain
$$a_0(n,\a)-a_0(n-1,\a)={n+\a \ch n-1},\;n=1,2,3,\ldots$$
and therefore, by using (\ref{binform})
$$a_0(n,\a)=a_0(n,\a)-a_0(0,\a)=\sum_{k=1}^n{k+\a \ch k-1}=
{n+\a+1 \ch n-1},\;n=1,2,3,\ldots,$$
which proves (\ref{anul}). Note that this proof of (\ref{anul}) is different
from the one given in \cite{Bav}.

In order to prove (\ref{ai}) we write instead of (\ref{a1}) and (\ref{a2})~:
\be\la{aa1}\sum_{i=1}^{\infty}a_i(x)D^i\L=
{n+\a \ch n}D^2\Ln-a_0(n,\a)\L,\;\ndots\ee
and
\be\la{aa2}\sum_{i=1}^{\infty}a_i(x)D^{i+1}\L=
-{n+\a \ch n-1}D^2\Ln-a_0(n,\a)D\L,\;\ndots.\ee
First we will prove that every solution of (\ref{aa2}) is also a solution of
(\ref{aa1}). Note that
\be\la{difa}a_0(n+1,\a)=a_0(n,\a)+{n+\a+1 \ch n},\;\ndots,\ee
which was already obtained before. Now suppose that
$\l\{a_i(x)\r\}_{i=1}^{\infty}$ is a solution of (\ref{aa2}). Then we
use (\ref{rel1}), (\ref{rel2}) and (\ref{difa}) to find for the left-hand
side of (\ref{aa1})
\bea & &\sum_{i=1}^{\infty}a_i(x)D^i\L\nn
&=&\sum_{i=1}^{\infty}a_i(x)D^{i+1}\L-
\sum_{i=1}^{\infty}a_i(x)D^{i+1}\Ln\nn
&=&{n+\a+1 \ch n}D^2L_{n+2}^{(\a)}(x)-{n+\a \ch n-1}D^2\Ln+{}\nn
& &{}\hspace{1cm}{}+a_0(n+1,\a)D\Ln-a_0(n,\a)D\L\nn
&=&{n+\a+1 \ch n}\l[D^2\Ln-D\Ln\r]-{n+\a \ch n-1}D^2\Ln\nn
& &{}\hspace{1cm}{}+\l[a_0(n,\a)+{n+\a+1 \ch n}\r]D\Ln-a_0(n,\a)D\L\nn
&=&\l[{n+\a+1 \ch n}-{n+\a \ch n-1}\r]D^2\Ln+a_0(n,\a)\l[D\Ln-D\L\r]\nn
&=&{n+\a \ch n}D^2\Ln-a_0(n,\a)\L,\;\ndots,\n\eea
which equals the right-hand side of (\ref{aa1}).

Now we will solve (\ref{aa2}). We have by using (\ref{anul}), (\ref{defhyp})
and (\ref{diff})
\bea & &-{n+\a \ch n-1}D^2\Ln-a_0(n,\a)D\L\nn
&=&{n+\a+1 \ch n-1}L_{n-1}^{(\a+1)}(x)-
{n+\a \ch n-1}L_{n-1}^{(\a+2)}(x)\nn
&=&-\frac{(n-1)x}{(\a+2)(\a+3)}{n+\a \ch n-1}{n+\a+1 \ch n-1}
\sum_{k=0}^{\infty}\frac{(-n+2)_k}{(\a+4)_k}\frac{x^k}{k!}\nn
&=&-\frac{x}{\a+2}{n+\a \ch n-1}L_{n-2}^{(\a+3)}(x),\;n=2,3,4,\ldots.\n\eea
Note that the right-hand side of (\ref{aa2}) equals zero for $n=0$ and
$n=1$, which is necessary for solvability. Hence, the system of equations
(\ref{aa2}) is equivalent to
$$\sum_{i=1}^{\infty}a_i(x)D^{i+1}\L=
-\frac{x}{\a+2}{n+\a \ch n-1}L_{n-2}^{(\a+3)}(x),\;n=2,3,4,\ldots.$$
Now we apply the inversion formula of lemma~5 to obtain
$$a_i(x)=\frac{(-1)^ix}{\a+2}\sum_{j=1}^i{j+\a+1 \ch j}
L_{i-j}^{(-\a-i-2)}(-x)L_{j-1}^{(\a+3)}(x),\;i=1,2,3,\ldots$$
and with (\ref{F}) and (\ref{lemF}) we obtain
\bea a_i(x)&=&(-1)^ixF_{i-1}(\a+3,\a+3;2;x)\nn
&=&\frac{1}{i!}\sum_{j=1}^i(-1)^{i+j+1}
{\a+1 \ch j-1}{\a+2 \ch i-j}(\a+3)_{i-j}x^j,\;i=1,2,3,\ldots,\n\eea
which equals (\ref{ai}).

\subsection{The computation of the coefficients
$\l\{b_i(x)\r\}_{i=0}^{\infty}$}

We use the systems of equations $S_5=0$ and $S_6=0$ for $\ndots$ to find the
coefficients $\l\{b_i(x)\r\}_{i=0}^{\infty}$. Since $\a>-1$ we have
${n+\a \ch n-1}\ne 0$ for $n=1,2,3,\ldots$. We use this to conclude that
these systems are equivalent to
\be\la{b1}\sum_{i=0}^{\infty}b_i(x)D^i\L=\frac{1}{\a+1}{n+\a \ch n-1}
\l[(n-1)D^2\Ln+2D^3\Ln\r]\ee
for $\ndots$ and
\bea\la{b2}& &(n-1)\sum_{i=0}^{\infty}b_i(x)D^{i+1}\L+
\sum_{i=0}^{\infty}b_i(x)D^{i+2}\L\nn
&=&-\frac{n(\a+2)-(\a+1)}{(\a+1)(\a+3)}{n+\a \ch n-2}
\l[(n-1)D^2\Ln+2D^3\Ln\r]\eea
for $n=1,2,3,\ldots$. For $n=0$ equation (\ref{b1}) leads to $b_0(0,\a)=0$
and for $n=1$ equation (\ref{b2}) is trivial. Further we conclude that
$b_0(1,\a)$ is arbitrary. Since $b_i(x)$ must be a polynomial in $x$ of
degree at most $i$ for each $i=1,2,3,\ldots$ we may write
\be\la{kopbdef}b_i(x)=k_i'x^i+\mbox{lower order terms},\;i=1,2,3,\ldots.\ee
Then we find by comparing the coefficients of highest degree in (\ref{b1})
and (\ref{b2})~:
$$\frac{b_0(n,\a)}{n!}+\sum_{i=1}^n\frac{k_i'}{(n-i)!}=0,\;n=1,2,3,\ldots$$
and
\be\la{kopbrel}\frac{b_0(n,\a)}{(n-1)!}+\sum_{i=1}^{n-1}\frac{k_i'}{(n-i-1)!}=
\frac{n(\a+2)-(\a+1)}{(\a+1)(\a+3)}{n+\a \ch n-2}\frac{1}{(n-1)!},
\;n=2,3,4,\ldots.\ee
Hence, since $k_i'$ is independent of $n$, we obtain
$$b_0(n,\a)-b_0(n-1,\a)=\frac{n(\a+2)-(\a+1)}{(\a+1)(\a+3)}{n+\a \ch n-2},
\;n=2,3,4,\ldots$$
and therefore, by using (\ref{binform})
\bea b_0(n,\a)-b_0(1,\a)&=&\sum_{k=2}^n\frac{k(\a+2)-(\a+1)}{(\a+1)(\a+3)}
{k+\a \ch k-2}\nn
&=&\frac{\a+2}{\a+1}\sum_{k=2}^n{k+\a \ch k-3}+
\frac{1}{\a+1}\sum_{k=2}^n{k+\a \ch k-2}\nn
&=&\frac{\a+2}{\a+1}{n+\a+1 \ch n-3}+\frac{1}{\a+1}{n+\a+1 \ch n-2}\nn
&=&\frac{n(\a+2)-\a}{(\a+1)(\a+4)}{n+\a+1 \ch n-2},\;n=1,2,3,\ldots.\n\eea
This proves that $b_0(n,\a)-b_0(1,\a)=\b_0(n,\a),\;n=1,2,3,\ldots$, where
$\b_0(n,\a)$ is given by (\ref{bnul}).

Now we will show that every solution of (\ref{b1}) also satisfies
(\ref{b2}). In order to do this we write (\ref{b1}) in the form
\be\la{bb1}\sum_{i=1}^{\infty}b_i(x)D^i\L=F_n(x),\;\ndots,\ee
where
$$F_n(x)=\frac{1}{\a+1}{n+\a \ch n-1}
\l[(n-1)D^2\Ln+2D^3\Ln\r]-b_0(n,\a)\L.$$
Suppose that $\l\{b_i(x)\r\}_{i=1}^{\infty}$ is a solution of
(\ref{bb1}). If we now write
$$\sum_{i=1}^{\infty}b_i(x)D^{i+1}\L=G_n(x),\;\ndots,$$
then we find by using (\ref{rel1}) that
$$G_n(x)-G_{n+1}(x)=\sum_{i=1}^{\infty}b_i(x)D^i\L=F_n(x),\;\ndots.$$
Hence, we have $G_0(x)=0$ and
$$G_n(x)=-\sum_{k=0}^{n-1}F_k(x),\;n=1,2,3,\ldots.$$
In a similar way we find
\bea\sum_{i=1}^{\infty}b_i(x)D^{i+2}\L
&=&-\sum_{m=0}^{n-1}G_m(x)=\sum_{m=1}^{n-1}\sum_{k=0}^{m-1}F_k(x)\nn
&=&\sum_{k=0}^{n-2}\sum_{m=k+1}^{n-1}F_k(x)=\sum_{k=0}^{n-2}(n-1-k)F_k(x)\nn
&=&(n-1)\sum_{k=0}^{n-1}F_k(x)-\sum_{k=0}^{n-1}kF_k(x),\;n=2,3,4,\ldots.
\n\eea
Hence
$$(n-1)\sum_{i=1}^{\infty}b_i(x)D^{i+1}\L+
\sum_{i=1}^{\infty}b_i(x)D^{i+2}\L=-\sum_{k=0}^{n-1}kF_k(x),
\;n=1,2,3,\ldots.$$
In view of (\ref{b2}) we have to prove that
\bea\la{IV}\sum_{k=0}^{n-1}kF_k(x)&=&
\frac{n(\a+2)-(\a+1)}{(\a+1)(\a+3)}{n+\a \ch n-2}
\l[(n-1)D^2\Ln+2D^3\Ln\r]+{}\nn
& &{}\hspace{1cm}{}+b_0(n,\a)\l[(n-1)D\L+D^2\L\r],\;n=1,2,3,\ldots.\eea
This is done by induction. For $n=1$ this formula is trivial. Suppose that
(\ref{IV}) holds for certain $n\in\{1,2,3,\ldots\}$. Then we have to show
that
\bea\sum_{k=0}^nkF_k(x)&=&\frac{n(\a+2)+1}{(\a+1)(\a+3)}{n+\a+1 \ch n-1}
\l[nD^2L_{n+2}^{(\a)}(x)+2D^3L_{n+2}^{(\a)}(x)\r]+{}\nn
& &{}\hspace{1cm}{}+b_0(n+1,\a)\l[nD\Ln+D^2\Ln\r].\n\eea
We find
\bea\sum_{k=0}^nkF_k(x)&=&nF_n(x)+\sum_{k=0}^{n-1}kF_k(x)\nn
&=&nF_n(x)+\frac{n(\a+2)-(\a+1)}{(\a+1)(\a+3)}{n+\a \ch n-2}
\l[(n-1)D^2\Ln+2D^3\Ln\r]+{}\nn
& &{}\hspace{1cm}{}+b_0(n,\a)\l[(n-1)D\L+D^2\L\r]\nn
&=&\frac{n}{\a+1}{n+\a \ch n-1}\l[(n-1)D^2\Ln+2D^3\Ln\r]-nb_0(n,\a)\L+{}\nn
& &{}+\frac{n(\a+2)-(\a+1)}{(\a+1)(\a+3)}{n+\a \ch n-2}
\l[(n-1)D^2\Ln+2D^3\Ln\r]+{}\nn
& &{}\hspace{1cm}{}+b_0(n,\a)\l[(n-1)D\L+D^2\L\r].\n\eea
Earlier we have found that
\be\la{difb}b_0(n+1,\a)=b_0(n,\a)+\frac{n(\a+2)+1}{(\a+1)(\a+3)}
{n+\a+1 \ch n-1},\;n=1,2,3,\ldots.\ee
We also have
$$\frac{n}{\a+1}{n+\a \ch n-1}+
\frac{n(\a+2)-(\a+1)}{(\a+1)(\a+3)}{n+\a \ch n-2}=
\frac{n(\a+2)+1}{(\a+1)(\a+3)}{n+\a+1 \ch n-1}.$$
Hence we find by using (\ref{rel2})
\bea\sum_{k=0}^nkF_k(x)
&=&\frac{n(\a+2)+1}{(\a+1)(\a+3)}{n+\a+1 \ch n-1}\times{}\nn
& &{}\hspace{1cm}{}\times\l[(n-1)D^2\Ln+2D^3\Ln-nD\Ln-D^2\Ln\r]+{}\nn
& &{}\hspace{4cm}{}+b_0(n+1,\a)\l[nD\Ln+D^2\Ln\r]\nn
&=&\frac{n(\a+2)+1}{(\a+1)(\a+3)}{n+\a+1 \ch n-1}
\l[nD^2L_{n+2}^{(\a)}(x)+2D^3L_{n+2}^{(\a)}(x)\r]+{}\nn
& &{}\hspace{4cm}{}+b_0(n+1,\a)\l[nD\Ln+D^2\Ln\r].\n\eea
This completes the proof.

Now we will solve (\ref{bb1}). Note that $F_0(x)=0$ which is necessary for
solvability. Substitution of
$b_0(n,\a)=b_0(1,\a)+\b_0(n,\a),\;n=1,2,3,\ldots$ into the right-hand side
of (\ref{bb1}) now leads to
$$F_n(x)=-b_0(1,\a)\L+H_n(x),\;n=1,2,3,\ldots,$$
where
$$H_n(x)=\frac{1}{\a+1}{n+\a \ch n-1}
\l[(n-1)D^2\Ln+2D^3\Ln\r]-\b_0(n,\a)\L,\;n=1,2,3,\ldots.$$
As before we have
\bea \b_0(n,\a)&=&\frac{n(\a+2)-\a}{(\a+1)(\a+4)}{n+\a+1 \ch n-2}\nn
&=&\frac{1}{\a+1}{n+\a+1 \ch n-2}+\frac{\a+2}{\a+1}{n+\a+1 \ch n-3},
\;n=1,2,3,\ldots.\n\eea
Hence, since $\b_0(1,\a)=0$, we find that $H_1(x)=0$. Further we find by
using (\ref{diff})
\bea H_n(x)&=&\frac{\a+2}{\a+1}{n+\a \ch n-2}L_{n-1}^{(\a+2)}(x)
-\frac{2}{\a+1}{n+\a \ch n-1}L_{n-2}^{(\a+3)}(x)+{}\nn
& &{}-\frac{1}{\a+1}{n+\a+1 \ch n-2}\L-\frac{\a+2}{\a+1}{n+\a+1 \ch n-3}\L,
\;n=2,3,4,\ldots.\n\eea
Applying the inversion formula of lemma~5 we obtain
$$b_i(x)=b_0(1,\a)b_i^*(x)+\b_i(x),\;i=1,2,3,\ldots,$$
where, by using (\ref{inv}) and (\ref{def2})
\bea b_i^*(x)&=&(-1)^{i+1}\sum_{j=1}^iL_{i-j}^{(-\a-i-1)}(-x)L_j^{(\a)}(x)
=(-1)^iL_i^{(-\a-i-1)}(-x)\nn
&=&\sum_{j=0}^i\frac{(\a+1)_{i-j}}{(i-j)!}\frac{(-x)^j}{j!}=
\frac{1}{i!}\sum_{j=0}^i(-1)^j{i \ch j}(\a+1)_{i-j}x^j,
\;i=1,2,3,\ldots,\n\eea
which equals (\ref{bster}). Further we obtain that $\b_1(x)=0$
and
$$\b_i(x):=\b_i^{(1)}(x)+\b_i^{(2)}(x)+\b_i^{(3)}(x)+
\b_i^{(4)}(x),\;i=2,3,4,\ldots,$$
where
$$\b_i^{(1)}(x)=(-1)^i\frac{\a+2}{\a+1}\sum_{j=2}^i{j+\a \ch j-2}
L_{i-j}^{(-\a-i-1)}(-x)L_{j-1}^{(\a+2)}(x),\;i=2,3,4,\ldots,$$
$$\b_i^{(2)}(x)=\frac{2(-1)^{i+1}}{\a+1}\sum_{j=2}^i{j+\a \ch j-1}
L_{i-j}^{(-\a-i-1)}(-x)L_{j-2}^{(\a+3)}(x),\;i=2,3,4,\ldots,$$
$$\b_i^{(3)}(x)=\frac{(-1)^{i+1}}{\a+1}\sum_{j=2}^i{j+\a+1 \ch j-2}
L_{i-j}^{(-\a-i-1)}(-x)L_j^{(\a)}(x),\;i=2,3,4,\ldots$$
and
$$\b_i^{(4)}(x)=(-1)^{i+1}\frac{\a+2}{\a+1}\sum_{j=2}^i{j+\a+1 \ch j-3}
L_{i-j}^{(-\a-i-1)}(-x)L_j^{(\a)}(x),\;i=2,3,4,\ldots.$$
Note that we have by using (\ref{F}) and (\ref{G})
$$\b_i^{(1)}(x)=\frac{(-1)^i}{\a+1}F_{i-1}(\a+2,\a+2;0;x),
\;i=2,3,4,\ldots,$$
$$\b_i^{(2)}(x)=(-1)^{i+1}\frac{2(\a+2)}{\a+1}F_{i-2}(\a+3,\a+3;2;x),
\;i=2,3,4,\ldots,$$
$$\b_i^{(3)}(x)=\frac{(-1)^{i+1}}{(\a+1)(\a+2)(\a+3)}G_i(\a+2,\a;-1;x),
\;i=2,3,4,\ldots$$
and
$$\b_i^{(4)}(x)=\frac{(-1)^{i+1}}{(\a+1)(\a+3)(\a+4)}G_i(\a+2,\a;-2;x),
\;i=2,3,4,\ldots.$$
Applying (\ref{lemF}) and (\ref{lemG}) we now easily obtain (\ref{bi1}),
(\ref{bi2}), (\ref{bi3}) and (\ref{bi4}). We remark that
$$x\b_i^{(2)}(x)=\frac{2(\a+2)}{\a+1}a_{i-1}(x),\;i=2,3,4,\ldots.$$

\subsection{The computation of the coefficients
$\l\{c_i(x)\r\}_{i=0}^{\infty}$}

We will use the systems of equations $S_3=0$, $S_4=0$, $S_7=0$ and $S_8=0$
for $\ndots$ to find the coefficients $\l\{c_i(x)\r\}_{i=0}^{\infty}$.
For $n=0$ these equations lead to $c_0(0,\a)=0$.

First of all we use (\ref{b2}) to write
\bea\la{bb2}& &\sum_{i=0}^{\infty}b_i(x)D^{i+2}\L=
-(n-1)\sum_{i=0}^{\infty}b_i(x)D^{i+1}\L+{}\nn
& &{}\hspace{1cm}{}-\frac{n(\a+2)-(\a+1)}{(\a+1)(\a+3)}{n+\a \ch n-2}
\l[(n-1)D^2\Ln+2D^3\Ln\r]\eea
for $\ndots$. Now we substitute (\ref{a1}), (\ref{a2}), (\ref{b1}) and
(\ref{bb2}) into the four systems of equations. Then we will show that
\bea\la{equiv}& &(n-1){n+\a+1 \ch n-1}\l[{n+\a+1 \ch n-1}S_3-nS_4\r]+{}\nn
& &\hspace{2cm}{}+(\a+1)(\a+3)\l[{n+\a+1 \ch n-1}S_7-nS_8\r]=0,\;\ndots.\eea
This can be done by straightforward but tedious computations as follows.
First we obtain
\bea & &{n+\a+1 \ch n-1}S_3-nS_4\nn
&=&{n+\a \ch n}{n+\a+1 \ch n-1}\sum_{i=0}^{\infty}b_i(x)D^{i+1}\L+{}\nn
& &{}-(\a+1){n+\a \ch n-2}\sum_{i=0}^{\infty}c_i(x)D^i\L
-n{n+\a \ch n}\sum_{i=0}^{\infty}c_i(x)D^{i+1}\L+{}\nn
& &{}+\frac{1}{\a+1}{n+\a \ch n-2}{n+\a \ch n}{n+\a+1 \ch n-1}\times{}\nn
& &{}\hspace{5cm}{}\times\l[(n-1)D^2\Ln+2D^3\Ln\r],\;\ndots\n\eea
and
\bea & &{n+\a+1 \ch n-1}S_7-nS_8\nn
&=&-\frac{1}{\a+1}{n+\a \ch n}{n+\a+1 \ch n-2}{n+\a+1 \ch n-1}
\sum_{i=0}^{\infty}b_i(x)D^{i+1}\L+{}\nn
& &{}+{n+\a \ch n-2}{n+\a+1 \ch n-2}\sum_{i=0}^{\infty}c_i(x)D^i\L+{}\nn
& &{}\hspace{3cm}{}+{n+\a \ch n-1}{n+\a+1 \ch n-2}
\sum_{i=0}^{\infty}c_i(x)D^{i+1}\L+{}\nn
& &{}-\frac{1}{(\a+1)^2}{n+\a \ch n-2}{n+\a \ch n}{n+\a+1 \ch n-2}
{n+\a+1 \ch n-1}\times{}\nn
& &{}\hspace{5cm}{}\times\l[(n-1)D^2\Ln+2D^3\Ln\r],\;\ndots,\n\eea
which eventually lead to (\ref{equiv}).

Now we use (\ref{equiv}) to conclude that for $n=1,2,3,\ldots$ the systems
of equations $S_3=0$, $S_4=0$, $S_7=0$ and $S_8=0$ are equivalent to
$$S_3=0,\;{n+\a+1 \ch n-1}S_3-nS_4=0\;\mbox{ and }\;S_8=0.$$
Finally this leads to the following systems of equations~:
\bea\la{c1}& &\sum_{i=0}^{\infty}c_i(x)D^i\L\nn
&=&-\frac{n\a-(\a+1)}{(\a+1)^2}{n+\a \ch n-2}{n+\a \ch n}D^2\Ln
-\frac{2}{\a+1}{n+\a \ch n-2}{n+\a \ch n}D^3\Ln+{}\nn
& &{}-\frac{1}{\a+1}{n+\a \ch n-1}\sum_{i=0}^{\infty}a_i(x)D^{i+2}\L
-{n+\a \ch n}\sum_{i=0}^{\infty}b_i(x)D^{i+1}\L,\eea
\bea\la{c2}& &\sum_{i=0}^{\infty}c_i(x)D^{i+1}\L\nn
&=&{n+\a \ch n-2}^2D^2\Ln+
\frac{2}{\a+1}{n+\a \ch n-2}{n+\a \ch n-1}D^3\Ln+{}\nn
& &{}+\frac{1}{\a+1}{n+\a \ch n-2}\sum_{i=0}^{\infty}a_i(x)D^{i+2}\L
+{n+\a \ch n-1}\sum_{i=0}^{\infty}b_i(x)D^{i+1}\L\eea
and
\bea\la{c3}& &\sum_{i=0}^{\infty}c_i(x)D^{i+2}\L\nn
&=&-\frac{n(\a+2)-(\a+1)}{(\a+1)(\a+3)}{n+\a \ch n-2}^2D^2\Ln
-\frac{2}{\a+1}{n+\a \ch n-2}^2D^3\Ln+{}\nn
& &{}-\frac{1}{\a+1}{n+\a \ch n-3}\sum_{i=0}^{\infty}a_i(x)D^{i+2}\L
-{n+\a \ch n-2}\sum_{i=0}^{\infty}b_i(x)D^{i+1}\L\eea
for $\ndots$. For $n=0$ we only find that $c_0(0,\a)=0$ and for $n=1$ we
find that $c_0(1,\a)=b_0(1,\a)$ and $c_1(x)=-b_0(1,\a)x$. Since $c_i(x)$
must be a polynomial in $x$ of degree at most $i$ for each $i=1,2,3,\ldots$
we may write
$$c_i(x)=k_i''x^i+\mbox{lower order terms},\;i=1,2,3,\ldots.$$
Then we find by comparing the coefficients of highest degree in (\ref{c1})
and (\ref{c2}) by using (\ref{kopbdef})
$$\frac{c_0(n,\a)}{n!}+\sum_{i=1}^n\frac{k_i''}{(n-i)!}=0,\;n=1,2,3,\ldots$$
and
\bea\frac{c_0(n,\a)}{(n-1)!}+\sum_{i=1}^{n-1}\frac{k_i''}{(n-i-1)!}&=&
-{n+\a \ch n-2}^2\frac{1}{(n-1)!}+{}\nn
& &{}+{n+\a \ch n-1}\l[\frac{b_0(n,\a)}{(n-1)!}+
\sum_{i=1}^{n-1}\frac{k_i'}{(n-i-1)!}\r],\;n=2,3,4,\ldots.\n\eea
Hence, since $k_i''$ is independent of $n$, we obtain by using
(\ref{kopbrel})
$$c_0(n,\a)-c_0(n-1,\a)=\frac{1}{(\a+1)(\a+2)}{n+\a \ch n-1}{n+\a+1 \ch n-2},
\;n=2,3,4,\ldots.$$
Since $c_0(1,\a)=b_0(1,\a)$ this proves that
$c_0(n,\a)=b_0(1,\a)+\c_0(n,\a),\;n=1,2,3,\ldots$, where
\be\la{cnulsom}\c_0(n,\a)=\frac{1}{(\a+1)(\a+2)}
\sum_{k=1}^n{k+\a \ch k-1}{k+\a+1 \ch k-2},\;n=1,2,3,\ldots.\ee
By using lemma~4 we see that this equals (\ref{cnul}).

Now we will show that every solution of (\ref{c3}) also satisfies (\ref{c2})
and that every solution of (\ref{c2}) also satisfies (\ref{c1}). In order to
do this we write (\ref{c1}), (\ref{c2}) and (\ref{c3}) in the form
$$\sum_{i=0}^{\infty}c_i(x)D^{i+k}\L=F_n^{(k)}(x),\;\ndots,\;k=0,1,2,$$
respectively, where
\bea (-1)^{k+1}F_n^{(k)}(x)&=&\frac{n(\a+k)-(\a+1)}{(\a+1)(\a+k+1)}
{n+\a \ch n-2}{n+\a \ch n-k}D^2\Ln+{}\nn
& &{}\hspace{1cm}{}+\frac{2}{\a+1}{n+\a \ch n-2}{n+\a \ch n-k}D^3\Ln+{}\nn
& &{}\hspace{1cm}{}+\frac{1}{\a+1}{n+\a \ch n-k-1}
\sum_{i=0}^{\infty}a_i(x)D^{i+2}\L+{}\nn
& &{}\hspace{1cm}{}+{n+\a \ch n-k}\sum_{i=0}^{\infty}b_i(x)D^{i+1}\L,
\;k=0,1,2.\n\eea
Suppose that $\l\{c_i(x)\r\}_{i=1}^{\infty}$ is a solution of (\ref{c3}).
Then we find from
\be\la{cc3}\sum_{i=1}^{\infty}c_i(x)D^{i+2}\L=
F_n^{(2)}(x)-c_0(n,\a)D^2\L,\;\ndots\ee
and (\ref{rel1}) that
\bea\sum_{i=1}^{\infty}c_i(x)D^{i+1}\L
&=&\sum_{i=1}^{\infty}c_i(x)D^{i+2}\L
-\sum_{i=1}^{\infty}c_i(x)D^{i+2}\Ln\nn
&=&F_n^{(2)}(x)-c_0(n,\a)D^2\L-F_{n+1}^{(2)}(x)+c_0(n+1,\a)D^2\Ln\n\eea
for $\ndots$. Now we use (\ref{rel2}) to obtain for $\ndots$
$$\sum_{i=0}^{\infty}c_i(x)D^{i+1}\L=F_n^{(2)}(x)-F_{n+1}^{(2)}(x)
+\l[c_0(n+1,\a)-c_0(n,\a)\r]D^2\Ln.$$
So it remains to show that the right-hand side equals $F_n^{(1)}(x)$. This
can be achieved by straightforward but tedious computations. First we use
(\ref{difa}), (\ref{rel2}) and (\ref{a2}) to find
\bea & &a_0(n+1,\a)D^2\Ln+\sum_{i=1}^{\infty}a_i(x)D^{i+2}\Ln\nn
&=&a_0(n,\a)\l[D^2\L-D\L\r]+{n+\a+1 \ch n}D^2\Ln+{}\nn
& &{}\hspace{5cm}{}+\sum_{i=1}^{\infty}a_i(x)\l[D^{i+2}\L-D^{i+1}\L\r]\nn
&=&\sum_{i=0}^{\infty}a_i(x)D^{i+2}\L-\sum_{i=0}^{\infty}a_i(x)D^{i+1}\L
+{n+\a+1 \ch n}D^2\Ln\nn
&=&\sum_{i=0}^{\infty}a_i(x)D^{i+2}\L+
\l[{n+\a \ch n-1}+{n+\a+1 \ch n}\r]D^2\Ln,\;\ndots.\n\eea
In a similar way we use (\ref{difb}), (\ref{rel2}) and (\ref{b1}) to find
\bea & &b_0(n+1,\a)D\Ln+\sum_{i=1}^{\infty}b_i(x)D^{i+1}\Ln\nn
&=&b_0(n,\a)\l[D\L-\L\r]
+\frac{n(\a+2)+1}{(\a+1)(\a+3)}{n+\a+1 \ch n-1}D\Ln+{}\nn
& &{}\hspace{5cm}{}+\sum_{i=1}^{\infty}b_i(x)\l[D^{i+1}\L-D^i\L\r]\nn
&=&\sum_{i=0}^{\infty}b_i(x)D^{i+1}\L-\sum_{i=0}^{\infty}b_i(x)D^i\L
+\frac{n(\a+2)+1}{(\a+1)(\a+3)}{n+\a+1 \ch n-1}D\Ln\nn
&=&\sum_{i=0}^{\infty}b_i(x)D^{i+1}\L
+\frac{n(\a+2)+1}{(\a+1)(\a+3)}{n+\a+1 \ch n-1}D\Ln+{}\nn
& &{}\hspace{2cm}{}-\frac{1}{\a+1}{n+\a \ch n-1}\l[(n-1)D^2\Ln+2D^3\Ln\r],
\;n=1,2,3,\ldots.\n\eea
Now we use (\ref{rel2}) to obtain
\bea & &F_n^{(2)}(x)-F_{n+1}^{(2)}(x)\nn
&=&\frac{1}{\a+1}{n+\a \ch n-2}\sum_{i=0}^{\infty}a_i(x)D^{i+2}\L+
{n+\a \ch n-1}\sum_{i=0}^{\infty}b_i(x)D^{i+1}\L+{}\nn
& &{}+\frac{2}{\a+1}{n+\a \ch n-2}{n+\a \ch n-1}D^3\Ln
+{n+\a \ch n-2}^2D^2\Ln+{}\nn
& &{}-\frac{1}{(\a+1)(\a+2)}{n+\a+1 \ch n}{n+\a+2 \ch n-1}D^2\Ln,
\;n=1,2,3,\ldots.\n\eea
Finally we use
\be\la{difc}c_0(n+1,\a)-c_0(n,\a)=\frac{1}{(\a+1)(\a+2)}
{n+\a+1 \ch n}{n+\a+2 \ch n-1},\;n=1,2,3,\ldots,\ee
which was already obtained before, to conclude that
$$F_n^{(2)}(x)-F_{n+1}^{(2)}(x)+\l[c_0(n+1,\a)-c_0(n,\a)\r]D^2\Ln
=F_n^{(1)}(x),\;n=1,2,3,\ldots.$$
For $n=0$ we easily find
$$F_0^{(2)}(x)-F_1^{(2)}(x)+\l[c_0(1,\a)-c_0(0,\a)\r]D^2L_1^{(\a)}(x)
=0=F_0^{(1)}(x).$$

Now suppose that $\l\{c_i(x)\r\}_{i=1}^{\infty}$ is a solution of
(\ref{c2}). Then we have by using (\ref{rel1}) for $\ndots$
\bea\sum_{i=1}^{\infty}c_i(x)D^i\L
&=&\sum_{i=1}^{\infty}c_i(x)D^{i+1}\L
-\sum_{i=1}^{\infty}c_i(x)D^{i+1}\Ln\nn
&=&F_n^{(1)}(x)-c_0(n,\a)D\L-F_{n+1}^{(1)}(x)+c_0(n+1,\a)D\Ln.\n\eea
We use (\ref{rel2}) again to obtain for $\ndots$
$$\sum_{i=0}^{\infty}c_i(x)D^i\L=F_n^{(1)}(x)-F_{n+1}^{(1)}(x)
+\l[c_0(n+1,\a)-c_0(n,\a)\r]D\Ln.$$
Now we have to show that the right-hand side equals $F_n^{(0)}(x)$. This can
also be achieved by straightforward but again tedious computations. As
before we obtain by using (\ref{rel2})
\bea & &F_n^{(1)}(x)-F_{n+1}^{(1)}(x)\nn
&=&-\frac{1}{\a+1}{n+\a \ch n-1}\sum_{i=0}^{\infty}a_i(x)D^{i+2}\L-
{n+\a \ch n}\sum_{i=0}^{\infty}b_i(x)D^{i+1}\L+{}\nn
& &{}-\frac{2}{\a+1}{n+\a \ch n-2}{n+\a \ch n}D^3\Ln
-\frac{n\a-(\a+1)}{(\a+1)^2}{n+\a \ch n-2}{n+\a \ch n}D^2\Ln+{}\nn
& &{}-\frac{1}{(\a+1)(\a+2)}{n+\a+1 \ch n}{n+\a+2 \ch n-1}D\Ln,
\;n=1,2,3,\ldots.\n\eea
Finally we use (\ref{difc}) again to conclude that
$$F_n^{(1)}(x)-F_{n+1}^{(1)}(x)+\l[c_0(n+1,\a)-c_0(n,\a)\r]D\Ln
=F_n^{(0)}(x),\;n=1,2,3,\ldots.$$
For $n=0$ we easily find
\bea & &F_0^{(1)}(x)-F_1^{(1)}(x)+\l[c_0(1,\a)-c_0(0,\a)\r]DL_1^{(\a)}(x)\nn
&=&-b_0(1,\a)DL_1^{(\a)}(x)+b_0(1,\a)DL_1^{(\a)}(x)=0=F_0^{(0)}(x).\n\eea

Now we will solve (\ref{cc3}). In order to reduce the number of terms
involved we use (\ref{b2}) again to find for $\ndots$
\bea & &{}-{n+\a \ch n-2}\sum_{i=0}^{\infty}b_i(x)D^{i+1}\L=
\frac{1}{\a+2}{n+\a \ch n-1}\sum_{i=0}^{\infty}b_i(x)D^{i+2}\L+{}\nn
& &{}+\frac{n(\a+2)-(\a+1)}{(\a+1)(\a+3)}{n+\a \ch n-2}
\l[{n+\a \ch n-2}D^2\Ln+\frac{2}{\a+2}{n+\a \ch n-1}D^3\Ln\r].\n\eea
We use this to write (\ref{cc3}) in the form
\bea & &\sum_{i=1}^{\infty}c_i(x)D^{i+2}\L\nn
&=&{}-c_0(n,\a)D^2\L-
\frac{2}{(\a+1)(\a+2)}{n+\a \ch n-3}{n+\a \ch n-1}D^3\Ln+{}\nn
& &{}-\frac{1}{\a+1}{n+\a \ch n-3}\sum_{i=0}^{\infty}a_i(x)D^{i+2}\L+
\frac{1}{\a+2}{n+\a \ch n-1}\sum_{i=0}^{\infty}b_i(x)D^{i+2}\L\n\eea
for $\ndots$. We remark that it is necessary for solvability that the
right-hand side equals zero for $n=0$, $n=1$ and $n=2$. For $n=0$ and $n=1$
this is trivial. For $n=2$ we use (\ref{diff}), (\ref{def1}),
(\ref{bsplit}), (\ref{csplit}), (\ref{bnul}) and (\ref{cnul}) to see that
this right-hand side equals
$$-c_0(2,\a)D^2L_2^{(\a)}(x)+b_0(2,\a)D^2L_2^{(\a)}(x)=
-\c_0(2,\a)+\b_0(2,\a)=0.$$
Now we substitute $c_0(n,\a)=b_0(1,\a)+\c_0(n,\a),\;n=3,4,5,\ldots$ to find
by using (\ref{bsplit})
\be\la{cinv}\sum_{i=1}^{\infty}c_i(x)D^{i+2}\L=b_0(1,\a)K_n(x)+M_n(x),
\;n=3,4,5,\ldots,\ee
where
$$K_n(x)=\l[\frac{1}{\a+2}{n+\a \ch n-1}-1\r]D^2\L
+\frac{1}{\a+2}{n+\a \ch n-1}\sum_{i=1}^{\infty}b_i^*(x)D^{i+2}\L$$
and
\bea & &M_n(x)
={}-\frac{2}{(\a+1)(\a+2)}{n+\a \ch n-3}{n+\a \ch n-1}D^3\Ln+{}\nn
& &{}\hspace{1cm}{}-\l[\c_0(n,\a)+\frac{1}{\a+1}{n+\a \ch n-3}a_0(n,\a)
-\frac{1}{\a+2}{n+\a \ch n-1}\b_0(n,\a)\r]D^2\L+{}\nn
& &{}\hspace{1cm}{}-\frac{1}{\a+1}{n+\a \ch n-3}
\sum_{i=1}^{\infty}a_i(x)D^{i+2}\L+
\frac{1}{\a+2}{n+\a \ch n-1}\sum_{i=1}^{\infty}\b_i(x)D^{i+2}\L.\n\eea

By using (\ref{bster}) we may write for $k\in\{0,1,2,\ldots\}$ and $\ndots$
$$\sum_{i=1}^{\infty}b_i^*(x)D^{i+k}\L=-D^k\L+\sum_{i=0}^{\infty}
\sum_{j=0}^i\frac{(-1)^j}{i!}{i \ch j}(\a+1)_{i-j}x^jD^{i+k}\L.$$
Changing the order of summation we find
\bea & &\sum_{i=0}^{\infty}\sum_{j=0}^i
\frac{(-1)^j}{i!}{i \ch j}(\a+1)_{i-j}x^jD^{i+k}\L\nn
&=&\sum_{j=0}^{\infty}\sum_{i=0}^{\infty}\frac{(-1)^j}{(i+j)!}{i+j \ch j}
(\a+1)_ix^jD^{i+j+k}\L\nn
&=&\sum_{i=0}^{\infty}\frac{(\a+1)_i}{i!}
\sum_{j=0}^{\infty}\frac{(-1)^j}{j!}x^jD^{i+j+k}\L,\;k=0,1,2,\ldots.\n\eea
Now we use (\ref{defhyp}) to obtain for $i=0,1,2,\ldots$
\bea\sum_{j=0}^{\infty}\frac{(-1)^j}{j!}x^jD^{i+j+k}\L
&=&{n+\a \ch n}\sum_{j=0}^{\infty}\frac{(-1)^j}{j!}x^j\sum_{m=i+j+k}^{\infty}
\frac{(-n)_m}{(\a+1)_m}\frac{x^{m-i-j-k}}{(m-i-j-k)!}\nn
&=&{n+\a \ch n}\sum_{j=0}^{\infty}\sum_{m=j}^{\infty}\frac{(-1)^j}{j!}
\frac{(-n)_{m+i+k}}{(\a+1)_{m+i+k}}\frac{x^m}{(m-j)!}\nn
&=&{n+\a \ch n}\sum_{m=0}^{\infty}\frac{(-n)_{m+i+k}}{(\a+1)_{m+i+k}}
\frac{x^m}{m!}\sum_{j=0}^m(-1)^j{m \ch j}\nn
&=&{n+\a \ch n}\frac{(-n)_{i+k}}{(\a+1)_{i+k}},\;k=0,1,2,\ldots.\n\eea
Hence, by using the Vandermonde summation formula (\ref{Van}) we find for
$n=1,2,3,\ldots$
\bea & &\sum_{i=0}^{\infty}\sum_{j=0}^i
\frac{(-1)^j}{i!}{i \ch j}(\a+1)_{i-j}x^jD^{i+k}\L\nn
&=&{n+\a \ch n}\sum_{i=0}^{\infty}
\frac{(\a+1)_i}{i!}\frac{(-n)_{i+k}}{(\a+1)_{i+k}}\nn
&=&{n+\a \ch n}\frac{(-n)_k}{(\a+1)_k}\hyp{2}{1}{-n+k,\a+1}{\a+k+1}{1}\nn
&=&{n+\a \ch n}\frac{(-n)_k}{(\a+1)_k}\frac{(k)_{n-k}}{(\a+k+1)_{n-k}}
=\frac{(-n)_k}{n\G(k)},\;k\in\{0,1,2,\ldots,n\}.\n\eea
This implies that for $n=1,2,3,\ldots$
$$\sum_{i=1}^{\infty}b_i^*(x)D^{i+k}\L=\frac{(-n)_k}{n\G(k)}-D^k\L,
\;k=0,1,2,\ldots.$$
Hence
$$\sum_{i=1}^{\infty}b_i^*(x)D^{i+2}\L=n-1-D^2\L,\;n=1,2,3,\ldots,$$
which implies that
$$K_n(x)={n+\a \ch n-2}-D^2\L,\;n=3,4,5,\ldots.$$

By using (\ref{anul}), (\ref{bnul}) and (\ref{cnulsom}) we obtain
\bea & &\c_0(n,\a)+\frac{1}{\a+1}{n+\a \ch n-3}a_0(n,\a)-
\frac{1}{\a+2}{n+\a \ch n-1}\b_0(n,\a)\nn
&=&\frac{1}{\a+1}\l[{n+\a \ch n-3}{n+\a+1 \ch n-1}-
\frac{n(\a+2)-\a}{(\a+2)(\a+4)}{n+\a \ch n-1}{n+\a+1 \ch n-2}\r]+{}\nn
& &{}\hspace{2cm}{}+\frac{1}{(\a+1)(\a+2)}
\sum_{k=1}^n{k+\a \ch k-1}{k+\a+1 \ch k-2}\nn
&=&\frac{2}{(\a+1)(\a+2)}{n+\a \ch n-1}{n+\a+1 \ch n-3}+{}\nn
& &{}\hspace{2cm}{}+\frac{1}{(\a+1)(\a+2)}
\sum_{k=1}^{n-1}{k+\a \ch k-1}{k+\a+1 \ch k-2},\;n=3,4,5,\ldots.\n\eea
Hence, we find by using (\ref{diff}) and (\ref{cnulsom})
\bea M_n(x)&=&{}-\frac{1}{\a+1}{n+\a \ch n-3}
\sum_{k=1}^{n-2}(-1)^ka_k(x)L_{n-k-2}^{(\a+k+2)}(x)+{}\nn
& &{}+\frac{1}{\a+2}{n+\a \ch n-1}
\sum_{k=1}^{n-2}(-1)^k\b_k(x)L_{n-k-2}^{(\a+k+2)}(x)+{}\nn
& &{}-\frac{2}{(\a+1)(\a+2)}{n+\a \ch n-1}
\l[{n+\a \ch n-3}D^3\Ln+{n+\a+1 \ch n-3}D^2\L\r]+{}\nn
& &{}-\c_0(n-1,\a)D^2\L,\;n=3,4,5,\ldots.\n\eea

Now we apply the inversion formula of lemma~5 to (\ref{cinv}) to find that
$$c_i(x)=b_0(1,\a)c_i^*(x)+\c_i(x),\;i=1,2,3,\ldots,$$
where for $i=1,2,3\ldots$ we have by using (\ref{diff})
\bea c_i^*(x)&=&(-1)^i\sum_{j=1}^iL_{i-j}^{(-\a-i-3)}(-x)
\l[{j+\a+2 \ch j}-L_j^{(\a+2)}(x)\r]\nn
&=&(-1)^i\l[\sum_{j=0}^i{j+\a+2 \ch j}L_{i-j}^{(-\a-i-3)}(-x)
-\sum_{j=0}^iL_{i-j}^{(-\a-i-3)}(-x)L_j^{(\a+2)}(x)\r].\n\eea
Now we use (\ref{rel}) and (\ref{def2}) to find for $i=1,2,3,\ldots$
$$\sum_{j=0}^i{j+\a+2 \ch j}L_{i-j}^{(-\a-i-3)}(-x)=
\sum_{j=0}^i(-1)^j{-\a-3 \ch j}L_{i-j}^{(-\a-i-3)}(-x)=
L_i^{(-i)}(-x)=\frac{x^i}{i!}.$$
This implies, by using (\ref{inv}), that
$$c_i^*(x)=(-1)^i\frac{x^i}{i!},\;i=1,2,3,\ldots,$$
which proves (\ref{cster}). Further we find
$$\c_i(x):=\c_i^{(1)}(x)+\c_i^{(2)}(x)+\c_i^{(3)}(x)
+\c_i^{(4)}(x)+\c_i^{(5)}(x),\;i=1,2,3,\ldots,$$
where for $i=1,2,3,\ldots$
$$\c_i^{(1)}(x)=\frac{(-1)^{i+1}}{\a+1}\sum_{k=1}^i(-1)^ka_k(x)
\sum_{j=k}^i{j+\a+2 \ch j-1}L_{i-j}^{(-\a-i-3)}(-x)L_{j-k}^{(\a+k+2)}(x),$$
$$\c_i^{(2)}(x)=\frac{(-1)^i}{\a+2}\sum_{k=1}^i(-1)^k\b_k(x)
\sum_{j=k}^i{j+\a+2 \ch j+1}L_{i-j}^{(-\a-i-3)}(-x)L_{j-k}^{(\a+k+2)}(x),$$
$$\c_i^{(3)}(x)=\frac{2(-1)^i}{(\a+1)(\a+2)}\sum_{j=1}^i{j+\a+2 \ch j-1}
{j+\a+2 \ch j+1}L_{i-j}^{(-\a-i-3)}(-x)L_j^{(\a+3)}(x),$$
$$\c_i^{(4)}(x)=\frac{2(-1)^{i+1}}{(\a+1)(\a+2)}\sum_{j=1}^i{j+\a+2 \ch j+1}
{j+\a+3 \ch j-1}L_{i-j}^{(-\a-i-3)}(-x)L_j^{(\a+2)}(x)$$
and
$$\c_i^{(5)}(x)=(-1)^{i+1}\sum_{j=1}^i\c_0(j+1,\a)
L_{i-j}^{(-\a-i-3)}(-x)L_j^{(\a+2)}(x).$$

Note that we have for $i=1,2,3,\ldots$ by using (\ref{F})
\be\la{ci1*}\c_i^{(1)}(x)=\frac{(-1)^{i+1}}{\a+1}
\sum_{k=1}^i(-1)^k(\a+4)_{k-1}a_k(x)F_{i-k}(\a+k+3,\a+k+2;k;x)\ee
and
\be\la{ci2*}\c_i^{(2)}(x)=(-1)^i\sum_{k=1}^i(-1)^k(\a+3)_k\b_k(x)
F_{i-k}(\a+k+3,\a+k+2;k+2;x).\ee
By using (\ref{lemF}) we easily obtain (\ref{ci1}) and (\ref{ci2}).

Further we have by using (\ref{H})
$$\c_i^{(3)}(x)=\frac{2(-1)^i}{(\a+1)(\a+3)}H_i(\a+3,\a+3,\a+3;0,2;x),
\;i=1,2,3,\ldots$$
and
$$\c_i^{(4)}(x)=\frac{2(-1)^{i+1}}{(\a+1)(\a+4)}H_i(\a+3,\a+4,\a+2;0,2;x),
\;i=1,2,3,\ldots.$$
Now we apply (\ref{lemH}) to obtain for $i=1,2,3,\ldots$
\bea\la{ci3*}& &\c_i^{(3)}(x)=\frac{2}{(\a+1)(\a+3)}\sum_{j=0}^i
\sum_{n=0}^j(-1)^j\frac{(-j)_n(\a+3)_{i-j+n}(\a+3)_n}{j!(i-j)!n!}\times{}\nn
& &{}\hspace{5cm}{}\times\hyphyp{3}{2}{-i+j,n+\a+3,n+\a+4}{n,n+2}{1}x^j\eea
and
\bea\la{ci4*}& &\c_i^{(4)}(x)=-\frac{2}{(\a+1)(\a+4)}\sum_{j=0}^i
\sum_{n=0}^j(-1)^j\frac{(-j)_n(\a+3)_{i-j+n}(\a+4)_n}{j!(i-j)!n!}\times{}\nn
& &{}\hspace{5cm}{}\times\hyphyp{3}{2}{-i+j,n+\a+4,n+\a+3}{n,n+2}{1}x^j,\eea
which lead to (\ref{ci3}) and (\ref{ci4}) by using (\ref{defphi}).

Finally we use (\ref{cnul}) to find for $i=1,2,3,\ldots$
\bea\c_i^{(5)}(x)&=&\frac{(-1)^{i+1}}{\a+1}\sum_{j=1}^i{j+\a+2 \ch j-1}
\hyp{3}{2}{-j+1,-\a-2,\a+3}{2,\a+4}{1}
L_{i-j}^{(-\a-i-3)}(-x)L_j^{(\a+2)}(x)\nn
&=&\frac{(-1)^{i+1}}{(\a+1)(\a+3)}\sum_{j=1}^i\sum_{k=0}^{j-1}
\frac{(-j+1)_k(-\a-2)_k(\a+3)_k}{(2)_k(\a+4)_kk!}\times{}\nn
& &{}\hspace{5cm}{}\times\frac{(\a+3)_j}{\G(j)}
L_{i-j}^{(-\a-i-3)}(-x)L_j^{(\a+2)}(x)\nn
&=&\frac{(-1)^{i+1}}{(\a+1)(\a+3)}\sum_{k=0}^{i-1}\sum_{j=k+1}^i
\frac{(-j+1)_k(-\a-2)_k(\a+3)_k}{(2)_k(\a+4)_kk!}\times{}\nn
& &{}\hspace{5cm}{}\times\frac{(\a+3)_j}{\G(j)}
L_{i-j}^{(-\a-i-3)}(-x)L_j^{(\a+2)}(x).\n\eea
Now we write
$$\frac{(-j+1)_k}{\G(j)}=(-1)^k\frac{(j-k)_k}{\G(j)}=\frac{(-1)^k}{\G(j-k)},
\;k=0,1,2,\ldots$$
to obtain by using (\ref{F}) and (\ref{lemF})
\bea\c_i^{(5)}(x)&=&\frac{(-1)^{i+1}}{(\a+1)(\a+3)}\sum_{k=0}^{i-1}(-1)^k
\frac{(-\a-2)_k(\a+3)_k}{(2)_k(\a+4)_kk!}\sum_{j=0}^i
\frac{(\a+3)_j}{\G(j-k)}L_{i-j}^{(-\a-i-3)}(-x)L_j^{(\a+2)}(x)\nn
&=&\frac{(-1)^{i+1}}{(\a+1)(\a+3)}\sum_{k=0}^{i-1}(-1)^k
\frac{(-\a-2)_k(\a+3)_k}{(2)_k(\a+4)_kk!}F_i(\a+3,\a+2;-k;x)\nn
&=&\frac{(-1)^{i+1}}{(\a+1)(\a+3)}\sum_{k=0}^{i-1}(-1)^k
\frac{(-\a-2)_k(\a+3)_k}{(2)_k(\a+4)_kk!}\times{}\nn
& &{}\hspace{3cm}{}\times\frac{1}{\G(i-k)}\sum_{j=0}^i
(-1)^j{k+\a+3 \ch j}{k+\a+3 \ch i-j}(\a+3)_{i-j}x^j,\n\eea
which proves (\ref{ci5}).

Now we will show that for nonnegative integer values of $\a$ we have
$$\c_i(x)=0\;\mbox{ for }\;i>4\a+10.$$

So we assume that $\a\in\{0,1,2,\ldots\}$. First we consider $\c_i^{(1)}(x)$
given by (\ref{ci1*}). We use (\ref{lemF}) to obtain
$$F_{i-k}(\a+k+3,\a+k+2;k;x)=\frac{1}{\G(i)}\sum_{j=0}^{i-k}
(-1)^j{\a+3 \ch j}{\a+3 \ch i-k-j}(\a+k+3)_{i-k-j}x^j.$$
The terms of this sum are equal to zero for $j>\a+3$ and for $i-k-j>\a+3$.
Hence all terms vanish if $i-k>2\a+6$. Since $a_k(x)=0$ for $k>2\a+4$ we
conclude from (\ref{ci1*}) in the same way that
$$\c_i^{(1)}(x)=0\;\mbox{ for }\;i>2\a+4+2\a+6=4\a+10.$$
Now we consider $\c_i^{(2)}(x)$ given by (\ref{ci2*}). In the same way we
have by using (\ref{lemF})
$$F_{i-k}(\a+k+3,\a+k+2;k+2;x)=\frac{1}{\G(i+2)}\sum_{j=0}^{i-k}
(-1)^j{\a+1 \ch j}{\a+1 \ch i-k-j}(\a+k+3)_{i-k-j}x^j.$$
The terms of this sum are equal to zero for $j>\a+1$ and for $i-k-j>\a+1$.
Hence all terms vanish if $i-k>2\a+2$. Since $\b_k(x)=0$ for $k>2\a+8$ we
conclude from (\ref{ci2*}) that
$$\c_i^{(2)}(x)=0\;\mbox{ for }\;i>2\a+8+2\a+2=4\a+10.$$

We need the following lemma~:

\lem{For $\ell\in\{0,1,2,\ldots\}$ and $b-a\notin\{1,2,3,\ldots\}$ we have
\be\la{lemphi}\hyphyp{3}{2}{-n,a,c+\ell}{b,c}{1}=\frac{1}{\G(n+b)}
\sum_{k=0}^{\ell}\frac{(-n)_k(-\ell)_k(a)_k(b-a-k)_n}{(a-b+1)_k\G(c+k)k!},
\;\ndots.\ee}

{\bf Proof.} Let $\ell\in\{0,1,2,\ldots\}$ and $b-a\notin\{1,2,3,\ldots\}$.
We have (see for instance \cite{Luke}, \S 9.1 formula (34))
$$\hyphyp{3}{2}{-n,a,c+\ell}{b,c}{1}=\sum_{k=0}^n(-1)^k
\frac{(-n)_k(-\ell)_k(a)_k}{\G(c+k)k!}\hyphyp{2}{1}{-n+k,a+k}{b+k}{1},
\;\ndots.$$
Hence, we obtain by using the Vandermonde summation formula (\ref{som})
\bea\hyphyp{3}{2}{-n,a,c+\ell}{b,c}{1}&=&\frac{1}{\G(n+b)}
\sum_{k=0}^n(-1)^k\frac{(-n)_k(-\ell)_k(a)_k(b-a)_{n-k}}{\G(c+k)k!}\nn
&=&\frac{1}{\G(n+b)}\sum_{k=0}^{\ell}
\frac{(-n)_k(-\ell)_k(a)_k(b-a-k)_n}{(a-b+1)_k\G(c+k)k!},\;\ndots\n\eea
which proves (\ref{lemphi}).

By using (\ref{ci3*}) and applying (\ref{lemphi}) twice we obtain for
$i=1,2,3,\ldots$
\bea\la{ci3**}\c_i^{(3)}(x)&=&\frac{2}{(\a+1)(\a+3)}
\sum_{k=0}^{\a+2}\frac{(-\a-2)_k(\a+3)_k}{(\a+4)_kk!}\times{}\nn
& &{}\times\sum_{j=0}^i(-1)^j
\frac{(-i+j)_k(\a+3)_{i-j}(-k-\a-3)_{i-j}}{j!(i-j)!}\times{}\nn
& &{}\hspace{5cm}{}\times\hyphyp{3}{2}{-j,i-j+\a+3,k+\a+3}{i-j,k+2}{1}x^j\nn
&=&\frac{2(-1)^i}{(\a+1)(\a+3)\G(i)}\sum_{k=0}^{\a+2}\sum_{m=0}^{\a+1}
\frac{(-\a-2)_k(\a+3)_k(-\a-1)_m(\a+3)_m}{(\a+4)_k(\a+4)_m\G(k+m+2)k!m!}
\times{}\nn
& &{}\times\sum_{j=0}^i(-1)^j{m+\a+3 \ch j}{k+\a+3 \ch i-j}
(-i+j)_k(-j)_m(m+\a+3)_{i-j}x^j.\eea
Note that the terms of the inner sum are equal to zero for $j>m+\a+3$ and
for $i-j>k+\a+3$. Hence all terms vanish if $i>k+m+2\a+6$. This implies
that $\c_i^{(3)}(x)=0$ if $i>k+m+2\a+6$ for all $k\in\{0,1,2,\ldots,\a+2\}$
and $m\in\{0,1,2,\ldots,\a+1\}$. Hence
$$\c_i^{(3)}(x)=0\;\mbox{ for }\;i>\a+2+\a+1+2\a+6=4\a+9.$$

In a similar way we use (\ref{ci4*}) and apply (\ref{lemphi}) twice again to
find for $i=1,2,3,\ldots$
\bea\la{ci4**}\c_i^{(4)}(x)&=&-\frac{2}{(\a+1)(\a+4)}
\sum_{k=0}^{\a+1}\frac{(-\a-1)_k(\a+4)_k}{(\a+5)_kk!}\times{}\nn
& &{}\times\sum_{j=0}^i(-1)^j
\frac{(-i+j)_k(\a+3)_{i-j}(-k-\a-4)_{i-j}}{j!(i-j)!}\times{}\nn
& &{}\hspace{5cm}{}\times\hyphyp{3}{2}{-j,i-j+\a+3,k+\a+4}{i-j,k+2}{1}x^j\nn
&=&\frac{2(-1)^{i+1}}{(\a+1)(\a+4)\G(i)}\sum_{k=0}^{\a+1}\sum_{m=0}^{\a+2}
\frac{(-\a-1)_k(\a+4)_k(-\a-2)_m(\a+3)_m}{(\a+5)_k(\a+4)_m\G(k+m+2)k!m!}
\times{}\nn
& &{}\times\sum_{j=0}^i(-1)^j{m+\a+3 \ch j}{k+\a+4 \ch i-j}
(-i+j)_k(-j)_m(m+\a+3)_{i-j}x^j.\eea
The terms of the inner sum are equal to zero for $j>m+\a+3$ and for
$i-j>k+\a+4$. Hence all terms vanish if $i>k+m+2\a+7$. This implies that
$$\c_i^{(4)}(x)=0\;\mbox{ for }\;i>\a+1+\a+2+2\a+7=4\a+10.$$

Further we note that for $i=1,2,3\ldots$ (\ref{ci5}) can now be written as
\bea\la{ci5*} & &\c_i^{(5)}(x)=\frac{(-1)^{i+1}}{(\a+1)(\a+3)}
\sum_{k=0}^{\a+2}(-1)^k\frac{(-\a-2)_k(\a+3)_k}{(2)_k(\a+4)_kk!}\times{}\nn
& &{}\hspace{3cm}{}\times\frac{1}{\G(i-k)}
\sum_{j=0}^i(-1)^j{k+\a+3 \ch j}{k+\a+3 \ch i-j}(\a+3)_{i-j}x^j.\eea
Now we see that the terms of the inner sum are equal to zero for
$j>k+\a+3$ and for $i-j>k+\a+3$. Hence all terms vanish if $i>2k+2\a+6$. So
we conclude that
$$\c_i^{(5)}(x)=0\;\mbox{ for }\;i>2\a+4+2\a+6=4\a+10.$$

Finally, it is not difficult to see that $\c_{4\a+10}^{(3)}(x)=0$ and that
$\c_{4\a+10}^{(p)}(x)$ reduces to one single term for $p\in\{1,2,4,5\}$.
Moreover, we find that
$$\c_{4\a+10}^{(1)}(x)+\c_{4\a+10}^{(2)}(x)=
\frac{-2x^{2\a+5}}{(\a+1)(\a+2)!(\a+4)!}=-\c_{4\a+10}^{(4)}(x)$$
and therefore
$$\c_{4\a+10}(x)=\c_{4\a+10}^{(5)}(x)=
\frac{x^{2\a+5}}{(\a+1)(2\a+5)(\a+2)!(\a+3)!}.$$

\section{The sums of the coefficients}

Earlier we have discovered (see \cite{Rap} and \cite{JCAM}) that
\be\la{soma}\sum_{i=1}^{\infty}a_i(x)=-\frac{\sin\pi\a}{\pi}
\frac{x}{(\a+2)(\a+3)}\hyp{1}{1}{1}{\a+4}{-x},\;\a>-1.\ee
In this section we will give a proof of (\ref{soma}). We will also prove
that
\be\la{somb}\sum_{i=1}^{\infty}\b_i(x)=\frac{\sin\pi\a}{\pi}
\frac{\a+2}{(\a+1)(\a+3)(\a+4)}\l[1-x\,\hyp{2}{2}{1,\a+3}{\a+2,\a+5}{-x}\r],
\;\a>-1.\ee

First we define for complex $x$
\be\la{K}K(p,q,r;s;\a,x):=\sum_{i=0}^{\infty}\sum_{j=0}^i
\frac{(-1)^{i+j}}{\G(i-s)}{\a+p \ch j}{\a+q \ch i-j}(\a+r)_{i-j}x^j,\ee
where $p,q,r$ and $s$ are integers with $p,q,r\ge 0$ and $\a>-1$. Then
we will prove the following lemma~:

\lem{For $\a\in\{0,1,2,\ldots\}$ we have
\be\la{lemK1}K(p,q,r;s;\a,x)=\sum_{j=0}^{\infty}{\a+p \ch j}
\frac{(j-\a-r-s)_{\a+q}}{\G(j+\a+q-s)}x^j\ee
for all complex $x$ and all integers $p,q,r$ and $s$ with $p,q,r\ge 0$.

For $\a>-1$ and $\a\notin\{0,1,2,\ldots\}$ the right-hand side of (\ref{K})
converges absolutely for all complex $x$ if $r+s<q$. In that case we have
\be\la{lemK2}K(p,q,r;s;\a,x)=\sum_{j=0}^{\infty}{\a+p \ch j}
\frac{\G(j+q-r-s)}{\G(j+\a+q-s)\G(j-\a-r-s)}x^j.\ee}

{\bf Proof.} First we assume that $\a\in\{0,1,2,\ldots\}$. Then the
right-hand side of (\ref{K}) is a finite sum since $\a+p$ and $\a+q$ are
nonnegative integers. So we may change the order of summation to obtain for
all complex $x$ and all integers $p,q,r$ and $s$ with $p,q,r\ge 0$~:
\bea\la{Kphi}K(p,q,r;s;\a,x)
&=&\sum_{i=0}^{\infty}\sum_{j=0}^i
\frac{(-1)^{i+j}}{\G(i-s)}{\a+p \ch j}{\a+q \ch i-j}(\a+r)_{i-j}x^j\nn
&=&\sum_{j=0}^{\infty}\sum_{i=0}^{\infty}
\frac{(-1)^i}{\G(i+j-s)}{\a+p \ch j}{\a+q \ch i}(\a+r)_ix^j\nn
&=&\sum_{j=0}^{\infty}{\a+p \ch j}x^j
\sum_{i=0}^{\infty}\frac{(-\a-q)_i(\a+r)_i}{\G(i+j-s)i!}\nn
&=&\sum_{j=0}^{\infty}{\a+p \ch j}\hyphyp{2}{1}{-\a-q,\a+r}{j-s}{1}x^j.\eea
Now we apply Vandermonde's summation formula (\ref{som}) to find
(\ref{lemK1}).

Now we assume that $\a>-1$ and $\a\notin\{0,1,2,\ldots\}$. Then (\ref{Kphi})
holds for all complex $x$ if $r+s<q$, since
$$\sum_{j=0}^{\infty}\sum_{i=0}^{\infty}{\a+p \ch j}
\frac{(-\a-q)_i(\a+r)_i}{\G(i+j-s)i!}x^j$$
converges absolutely for all complex $x$ if $r+s<q$.

In order to prove this, let $x$ be complex and define
$$u_{ij}:={\a+p \ch j}\frac{\G(i-\a-q)\G(i+\a+r)}
{\G(-\a-q)\G(\a+r)\G(i+j-s)i!}x^j,\;i,j=0,1,2,\ldots.$$
Now we will show that
$$\sum_{j=0}^{\infty}\sum_{i=0}^{\infty}\l|u_{ij}\r|$$
converges if $r+s<q$.

Since
$$(-1)^j{\a+p \ch j}=\frac{(-\a-p)_j}{j!}=
\frac{\G(j-\a-p)}{\G(-\a-p)j!}\sim\frac{j^{-\a-p-1}}{\G(-\a-p)}
\;\mbox{ for }\;j\rightarrow\infty$$
and $\a+p+1>p\ge 0$ we have
$$\lim\limits_{j\rightarrow\infty}{\a+p \ch j}=0.$$
Hence there exists a positive number $M_{\a,p}$ independent of $j$ such that
$$\l|{\a+p \ch j}\r|\le M_{\a,p},\;j=0,1,2,\ldots.$$
This implies that
\be\la{afs}\l|u_{ij}\r|\le C_{\a,p,q,r}
\l|\frac{\G(i-\a-q)\G(i+\a+r)}{\G(i+j-s)i!}\r||x|^j,\;i,j=0,1,2,\ldots,\ee
where
$$C_{\a,p,q,r}:=\frac{M_{\a,p}}{\l|\G(-\a-q)\G(\a+r)\r|}>0.$$

We distinguish two cases~: $s<0$ and $s\ge 0$.

If $s<0$ we have $\G(i+j-s)=(i+j-s-1)!\ge j!(i-s-1)!$. So we obtain
$$\l|u_{ij}\r|\le C_{\a,p,q,r}
\l|\frac{\G(i-\a-q)\G(i+\a+r)}{\G(i-s)\G(i+1)}\r|\frac{|x|^j}{j!},
\;i,j=0,1,2,\ldots.$$
Since
$$\l|\frac{\G(i-\a-q)\G(i+\a+r)}{\G(i-s)\G(i+1)}\r|\sim i^{r+s-q-1}
\;\mbox{ for }\;i\rightarrow\infty$$
we conclude that
$$\sum_{j=0}^{\infty}\sum_{i=0}^{\infty}\l|u_{ij}\r|<\infty$$
for $s<0$ if $r+s<q$.

If $s\ge 0$ we find by using (\ref{afs}) that $\l|u_{ij}\r|=0$ if
$i+j\le s$. So we may write
$$\sum_{j=0}^{\infty}\sum_{i=0}^{\infty}\l|u_{ij}\r|=
\sum_{j=0}^s\sum_{i=s+1-j}^{\infty}\l|u_{ij}\r|+
\sum_{j=s+1}^{\infty}\sum_{i=0}^{\infty}\l|u_{ij}\r|.$$

For $0\le j\le s$ we have
$$\l|\frac{\G(i-\a-q)\G(i+\a+r)}{\G(i+j-s)i!}\r|\sim i^{r+s-q-1-j}
\;\mbox{ for }\;i\rightarrow\infty$$
and $r+s-q-1-j\le r+s-q-1$. Hence
$$\sum_{j=0}^s\sum_{i=s+1-j}^{\infty}\l|u_{ij}\r|<\infty$$
if $r+s<q$.

For $j\ge s+1$ we may write
$\G(i+j-s)=(i+j-s-1)!\ge i!(j-s-1)!$. Hence
$$\l|u_{ij}\r|\le C_{\a,p,q,r}
\l|\frac{\G(i-\a-q)\G(i+\a+r)}{i!\,i!}\r|\frac{|x|^j}{(j-s-1)!}$$
for $i=0,1,2,\ldots$ and $j=s+1,s+2,s+3,\ldots$. Since
$$\l|\frac{\G(i-\a-q)\G(i+\a+r)}{i!\,i!}\r|\sim i^{r-q-2}
\;\mbox{ for }\;i\rightarrow\infty$$
and $r-q-2=r+s-q-1-s-1<r+s-q-1$ we conclude that we also have
$$\sum_{j=s+1}^{\infty}\sum_{i=0}^{\infty}\l|u_{ij}\r|<\infty$$
if $r+s<q$.

This implies that
$$\sum_{j=0}^{\infty}\sum_{i=0}^{\infty}\l|u_{ij}\r|<\infty$$
for $s\ge 0$ if $r+s<q$.

This proves that if $\a>-1$ and $\a\notin\{0,1,2,\ldots\}$
$$\sum_{j=0}^{\infty}\sum_{i=0}^{\infty}\l|u_{ij}\r|$$
converges if $r+s<q$. So (\ref{Kphi}) holds for all complex $x$ if
$r+s<q$. Now we apply Gauss' summation formula
$$\hyphyp{2}{1}{a,b}{c}{1}=\frac{\G(c-a-b)}{\G(c-a)\G(c-b)},\;c-a-b>0$$
to find (\ref{lemK2}). This proves lemma~9.

Note that (\ref{lemK2}) also holds for $\a\in\{0,1,2,\ldots\}$ if $r+s<q$.
So for $p=r+s\ge 0$ and $q=r+s+1$ we use (\ref{lemK2}) and the well-known
formula
\be\la{sin}\frac{1}{\G(z)\G(1-z)}=\frac{\sin\pi z}{\pi}\ee
to find for all $\a>-1$
$$K(r+s,r+s+1,r;s;\a,x)=\sum_{j=0}^{\infty}
\frac{\sin\pi(j-\a-r-s)}{\pi}\frac{\G(\a+r+s+1)}{\G(j+\a+r+1)}x^j.$$
Now we use
$$\sin\pi(j-\a-r-s)=(-1)^{j+r+s+1}\sin\pi\a,\;j=0,1,2,\ldots$$
to obtain for all $\a>-1$
\be\la{Kspec2}K(r+s,r+s+1,r;s;\a,x)=(-1)^{r+s+1}\frac{\sin\pi\a}{\pi}
\sum_{j=0}^{\infty}\frac{\G(\a+r+s+1)}{\G(j+\a+r+1)}(-x)^j\ee
for all integers $r$ and $s$ with $r\ge 0$ and $r+s\ge 0$.

Now we find from (\ref{ai}) for all $\a>-1$, by using (\ref{K}) and
(\ref{Kspec2})
\bea\sum_{i=1}^{\infty}a_i(x)&=&-xK(1,2,3;-2;\a,x)=-x\frac{\sin\pi\a}{\pi}
\sum_{j=0}^{\infty}\frac{\G(\a+2)}{\G(j+\a+4)}(-x)^j\nn
&=&-\frac{\sin\pi\a}{\pi}\frac{x}{(\a+2)(\a+3)}\hyp{1}{1}{1}{\a+4}{-x}.\n\eea
This proves (\ref{soma}).

In a similar way we obtain from (\ref{bi1}) and (\ref{bi2}) for all $\a>-1$,
by using (\ref{K}) and (\ref{Kspec2})
\bea\sum_{i=2}^{\infty}\b_i^{(1)}(x)&=&-\frac{1}{\a+1}K(2,3,2;0;\a,x)
=\frac{1}{\a+1}\frac{\sin\pi\a}{\pi}
\sum_{j=0}^{\infty}\frac{\G(\a+3)}{\G(j+\a+3)}(-x)^j\nn
&=&\frac{\sin\pi\a}{\pi}\frac{1}{\a+1}\hyp{1}{1}{1}{\a+3}{-x}\n\eea
and
\bea\sum_{i=2}^{\infty}\b_i^{(2)}(x)&=&-\frac{2(\a+2)}{\a+1}K(1,2,3;-2;\a,x)
=-\frac{2(\a+2)}{\a+1}\frac{\sin\pi\a}{\pi}\sum_{j=0}^{\infty}
\frac{\G(\a+2)}{\G(j+\a+4)}(-x)^j\nn
&=&-\frac{\sin\pi\a}{\pi}\frac{2}{(\a+1)(\a+3)}\hyp{1}{1}{1}{\a+4}{-x}.\n\eea

If $p=q=r+s+1$ we also obtain from (\ref{lemK2}) for all $\a>-1$
\be\la{Kspec3}K(r+s+1,r+s+1,r;s;\a,x)=\sum_{j=0}^{\infty}
\frac{(-\a-r-s-1)_j}{\G(j+\a+r+1)\G(j-\a-r-s)}(-x)^j\ee
for all integers $r$ and $s$ with $r\ge 0$ and $r+s+1\ge 0$.

Now we obtain from (\ref{bi3}) and (\ref{bi4}), by using (\ref{K}),
(\ref{Kspec3}) and (\ref{sin}) for all $\a>-1$
\bea\sum_{i=2}^{\infty}\b_i^{(3)}(x)&=&-\frac{1}{(\a+1)(\a+2)(\a+3)}
K(3,3,1;1;\a,x)+\frac{x}{(\a+1)(\a+2)}K(2,2,2;-1;\a,x)\nn
&=&-\frac{1}{(\a+1)(\a+2)(\a+3)}\sum_{j=0}^{\infty}
\frac{(-\a-3)_j}{\G(j+\a+2)\G(j-\a-2)}(-x)^j+{}\nn
& &\hspace{2cm}{}+\frac{x}{(\a+1)(\a+2)}\sum_{j=0}^{\infty}
\frac{(-\a-2)_j}{\G(j+\a+3)\G(j-\a-1)}(-x)^j\nn
&=&-\frac{1}{(\a+1)(\a+2)(\a+3)}\sum_{j=0}^{\infty}
\frac{(-\a-3)_j}{\G(j+\a+2)\G(j-\a-2)}(-x)^j+{}\nn
& &\hspace{2cm}{}+\frac{1}{(\a+1)(\a+2)(\a+3)}\sum_{j=1}^{\infty}
\frac{(-\a-3)_j}{\G(j+\a+2)\G(j-\a-2)}(-x)^j\nn
&=&-\frac{1}{(\a+1)(\a+2)(\a+3)}\frac{1}{\G(\a+2)\G(-\a-2)}
=\frac{\sin\pi\a}{\pi}\frac{1}{(\a+1)(\a+3)}\n\eea
and
\bea\sum_{i=2}^{\infty}\b_i^{(4)}(x)&=&-\frac{1}{(\a+1)(\a+3)(\a+4)}
K(4,4,1;2;\a,x)+\frac{x}{(\a+1)(\a+3)}K(3,3,2;0;\a,x)\nn
&=&-\frac{1}{(\a+1)(\a+3)(\a+4)}\sum_{j=0}^{\infty}
\frac{(-\a-4)_j}{\G(j+\a+2)\G(j-\a-3)}(-x)^j+{}\nn
& &\hspace{2cm}{}+\frac{x}{(\a+1)(\a+3)}\sum_{j=0}^{\infty}
\frac{(-\a-3)_j}{\G(j+\a+3)\G(j-\a-2)}(-x)^j\nn
&=&-\frac{1}{(\a+1)(\a+3)(\a+4)}\sum_{j=0}^{\infty}
\frac{(-\a-4)_j}{\G(j+\a+2)\G(j-\a-3)}(-x)^j+{}\nn
& &\hspace{2cm}{}+\frac{1}{(\a+1)(\a+3)(\a+4)}\sum_{j=1}^{\infty}
\frac{(-\a-4)_j}{\G(j+\a+2)\G(j-\a-3)}(-x)^j\nn
&=&-\frac{1}{(\a+1)(\a+3)(\a+4)}\frac{1}{\G(\a+2)\G(-\a-3)}
=-\frac{\sin\pi\a}{\pi}\frac{\a+2}{(\a+1)(\a+4)}.\n\eea

Hence, since $\b_1(x)=0$, we have
\bea & &\sum_{i=1}^{\infty}\b_i(x)=\sum_{i=2}^{\infty}
\l[\b_i^{(1)}(x)+\b_i^{(2)}(x)+\b_i^{(3)}(x)+\b_i^{(4)}(x)\r]\nn
&=&\frac{\sin\pi\a}{\pi}\frac{1}{\a+1}\l[\hyp{1}{1}{1}{\a+3}{-x}
-\frac{2}{\a+3}\hyp{1}{1}{1}{\a+4}{-x}-\frac{\a^2+4\a+2}{(\a+3)(\a+4)}\r]\nn
&=&\frac{\sin\pi\a}{\pi}\frac{1}{(\a+1)(\a+3)}\l[\frac{\a+2}{\a+4}+
\sum_{k=1}^{\infty}\frac{\a+k+1}{(\a+4)_k}(-x)^k\r]\nn
&=&\frac{\sin\pi\a}{\pi}\frac{\a+2}{(\a+1)(\a+3)(\a+4)}
\l[1-x\,\hyp{2}{2}{1,\a+3}{\a+2,\a+5}{-x}\r],\;\a>-1,\n\eea
which proves (\ref{somb}).

Note that for nonnegative integer values of $\a$ this implies that
$$\sum_{i=1}^{2\a+4}a_i(x)=0\;\mbox{ and }\;\sum_{i=1}^{2\a+8}\b_i(x)=0.$$

Finally we will show that for nonnegative integer values of $\a$ we also
have
$$\sum_{i=1}^{4\a+10}\c_i(x)=0.$$

Therefore we assume that $\a\in\{0,1,2,\ldots\}$. Now we will use
(\ref{lemK1}). Note that if $p=q=r+s\ge 0$ we obtain
\bea\la{Kspec4}K(r+s,r+s,r;s;\a,x)&=&\sum_{j=0}^{\infty}
\frac{(-\a-r-s)_{j+\a+r+s}}{\G(j+\a+r)j!}(-x)^j=
\frac{(-\a-r-s)_{\a+r+s}}{\G(\a+r)}\nn
&=&(-1)^{\a+r+s}\frac{\G(\a+r+s+1)}{\G(\a+r)}
=(-1)^{\a+r+s}\frac{(\a+1)_{r+s}}{(\a+1)_{r-1}}\eea
for all integers $r$ and $s$ with $r\ge 1$ and $r+s\ge 0$.

From (\ref{ci1*}) we easily obtain by changing the order of summation
$$\sum_{i=1}^{\infty}\c_i^{(1)}(x)=\frac{1}{\a+1}\sum_{k=1}^{\infty}
(\a+4)_{k-1}a_k(x)\sum_{i=0}^{\infty}(-1)^{i+1}F_i(\a+k+3,\a+k+2;k;x).$$
Note that all summations are in fact finite. Now we use (\ref{lemF}),
(\ref{K}) and (\ref{Kspec4}) to obtain
\bea & &\sum_{i=0}^{\infty}(-1)^iF_i(\a+k+3,\a+k+2;k;x)\nn
&=&\sum_{i=0}^{\infty}\sum_{j=0}^i\frac{(-1)^{i+j}}{\G(k+i)}
{\a+3 \ch j}{\a+3 \ch i-j}(\a+k+3)_{i-j}x^j\nn
&=&K(3,3,k+3;-k;\a,x)=\frac{(-1)^{\a+1}}{(\a+4)_{k-1}},
\;k=1,2,3,\ldots.\n\eea
Hence
$$\sum_{i=1}^{\infty}\c_i^{(1)}(x)=
\frac{(-1)^{\a}}{\a+1}\sum_{k=1}^{\infty}a_k(x)=0.$$
In the same way we find from (\ref{ci2*}) by changing the order of summation
$$\sum_{i=1}^{\infty}\c_i^{(2)}(x)=\sum_{k=1}^{\infty}
(\a+3)_k\b_k(x)\sum_{i=0}^{\infty}(-1)^iF_i(\a+k+3,\a+k+2;k+2;x).$$
Again we use (\ref{lemF}), (\ref{K}) and (\ref{Kspec4}) to find
\bea & &\sum_{i=0}^{\infty}(-1)^iF_i(\a+k+3,\a+k+2;k+2;x)\nn
&=&\sum_{i=0}^{\infty}\sum_{j=0}^i\frac{(-1)^{i+j}}{\G(k+i+2)}
{\a+1 \ch j}{\a+1 \ch i-j}(\a+k+3)_{i-j}x^j\nn
&=&K(1,1,k+3;-k-2;\a,x)=\frac{(-1)^{\a+1}}{(\a+2)_{k+1}},
\;k=1,2,3,\ldots.\n\eea
Hence
$$\sum_{i=1}^{\infty}\c_i^{(2)}(x)=
\frac{(-1)^{\a+1}}{\a+2}\sum_{k=1}^{\infty}\b_k(x)=0.$$

By using (\ref{ci3**}) we obtain
\bea & &\sum_{i=1}^{\infty}\c_i^{(3)}(x)\nn
&=&\frac{2}{(\a+1)(\a+3)}\sum_{k=0}^{\a+2}\sum_{m=0}^{\a+1}
\frac{(-\a-2)_k(\a+3)_k(-\a-1)_m(\a+3)_m}
{(\a+4)_k(\a+4)_m\G(k+m+2)k!m!}\times{}\nn
& &{}\hspace{1cm}{}\times\sum_{i=0}^{\infty}\sum_{j=0}^i
\frac{(-1)^{i+j}}{\G(i)}
{m+\a+3 \ch j}{k+\a+3 \ch i-j}(-i+j)_k(-j)_m(m+\a+3)_{i-j}x^j.\n\eea

Changing the order of summation we find
\bea & &\sum_{i=0}^{\infty}\sum_{j=0}^i\frac{(-1)^{i+j}}{\G(i)}
{m+\a+3 \ch j}{k+\a+3 \ch i-j}(-i+j)_k(-j)_m(m+\a+3)_{i-j}x^j\nn
&=&\sum_{j=0}^{\infty}{m+\a+3 \ch j}(-j)_mx^j\sum_{i=0}^{\infty}
\frac{(-1)^i}{\G(i+j)}{k+\a+3 \ch i}(-i)_k(m+\a+3)_i.\n\eea
Since $(-i)_k=0$ for $i<k$ we have by using the Vandermonde summation
formula (\ref{som})
\bea\la{S}S(j,k,m;\a)&:=&\sum_{i=0}^{\infty}
\frac{(-1)^i}{\G(i+j)}{k+\a+3 \ch i}(-i)_k(m+\a+3)_i\nn
&=&\sum_{i=k}^{\infty}\frac{(-k-\a-3)_i(-i)_k(m+\a+3)_i}{\G(i+j)i!}\nn
&=&(-1)^k(-k-\a-3)_k(m+\a+3)_k\,\hyphyp{2}{1}{-\a-3,k+m+\a+3}{k+j}{1}\nn
&=&(\a+4)_k(m+\a+3)_k\frac{(j-m-\a-3)_{\a+3}}{\G(k+j+\a+3)},
\;k=0,1,2,\ldots\eea
for every arbitrary $j$ and $m$ and for each $\a\in\{0,1,2,\ldots\}$. Hence
\bea & &\sum_{i=0}^{\infty}\sum_{j=0}^i\frac{(-1)^{i+j}}{\G(i)}
{m+\a+3 \ch j}{k+\a+3 \ch i-j}(-i+j)_k(-j)_m(m+\a+3)_{i-j}x^j\nn
&=&(\a+4)_k(m+\a+3)_k\sum_{j=0}^{\infty}{m+\a+3 \ch j}
\frac{(-j)_m(j-m-\a-3)_{\a+3}}{\G(k+j+\a+3)}x^j\nn
&=&(\a+4)_k(m+\a+3)_k\sum_{j=0}^{\infty}(-1)^j
\frac{(-j)_m(-m-\a-3)_{j+\a+3}}{\G(k+j+\a+3)j!}x^j\nn
&=&(\a+4)_k(m+\a+3)_k(-1)^m
\frac{(-m)_m(-m-\a-3)_{m+\a+3}}{\G(k+m+\a+3)m!}x^m\nn
&=&(\a+4)_k(m+\a+3)(-1)^{m+\a+1}x^m,\;k,m=0,1,2,\ldots,\n\eea
since $(-j)_m(-m-\a-3)_{j+\a+3}=0$ for all $j$ except $j=m$.

This implies that
\bea\sum_{i=1}^{\infty}\c_i^{(3)}(x)&=&\frac{2(-1)^{\a+1}}{(\a+1)(\a+3)}
\sum_{k=0}^{\a+2}\sum_{m=0}^{\a+1}(-1)^m\times{}\nn
& &{}\hspace{1cm}{}\times
\frac{(-\a-2)_k(\a+3)_k(-\a-1)_m(\a+3)_m}{(\a+4)_m\G(k+m+2)k!m!}
(m+\a+3)x^m\nn
&=&\frac{2(-1)^{\a+1}}{\a+1}\sum_{m=0}^{\a+1}(-1)^m\frac{(-\a-1)_m}{m!}
x^m\sum_{k=0}^{\a+2}\frac{(-\a-2)_k(\a+3)_k}{\G(k+m+2)k!}.\n\eea
Finally we use the Vandermonde summation formula (\ref{som}) again to obtain
$$\sum_{k=0}^{\a+2}\frac{(-\a-2)_k(\a+3)_k}{\G(k+m+2)k!}=
\hyphyp{2}{1}{-\a-2,\a+3}{m+2}{1}=\frac{(m-\a-1)_{\a+2}}{\G(m+\a+4)}.$$
Hence
$$\sum_{i=1}^{\infty}\c_i^{(3)}(x)=\frac{2(-1)^{\a+1}}{\a+1}
\sum_{m=0}^{\a+1}\frac{(-\a-1)_{m+\a+2}}{\G(m+\a+4)m!}(-x)^m=0.$$

Furthermore, by using (\ref{ci4**}) we obtain
\bea & &\sum_{i=1}^{\infty}\c_i^{(4)}(x)\nn
&=&\frac{2}{(\a+1)(\a+4)}\sum_{k=0}^{\a+1}\sum_{m=0}^{\a+2}
\frac{(-\a-1)_k(\a+4)_k(-\a-2)_m(\a+3)_m}
{(\a+5)_k(\a+4)_m\G(k+m+2)k!m!}\times{}\nn
& &{}\hspace{1cm}{}\times\sum_{i=0}^{\infty}\sum_{j=0}^i
\frac{(-1)^{i+j+1}}{\G(i)}
{m+\a+3 \ch j}{k+\a+4 \ch i-j}(-i+j)_k(-j)_m(m+\a+3)_{i-j}x^j.\n\eea

Changing the order of summation we find by using (\ref{S})
\bea & &\sum_{i=0}^{\infty}\sum_{j=0}^i\frac{(-1)^{i+j+1}}{\G(i)}
{m+\a+3 \ch j}{k+\a+4 \ch i-j}(-i+j)_k(-j)_m(m+\a+3)_{i-j}x^j\nn
&=&-\sum_{j=0}^{\infty}{m+\a+3 \ch j}(-j)_mx^j\sum_{i=0}^{\infty}
\frac{(-1)^i}{\G(i+j)}{k+\a+4 \ch i}(-i)_k(m+\a+3)_i\nn
&=&-\sum_{j=0}^{\infty}{m+\a+3 \ch j}(-j)_mx^jS(j,k,m-1;\a+1)\nn
&=&-(\a+5)_k(m+\a+3)_k\sum_{j=0}^{\infty}{m+\a+3 \ch j}
\frac{(-j)_m(j-m-\a-3)_{\a+4}}{\G(k+j+\a+4)}x^j\nn
&=&-(\a+5)_k(m+\a+3)_k\sum_{j=0}^{\infty}(-1)^j
\frac{(-j)_m(-m-\a-3)_{j+\a+4}}{\G(k+j+\a+4)j!}x^j,\;k,m=0,1,2,\ldots.\n\eea
Since $(-j)_m(-m-\a-3)_{j+\a+4}=0$ for all $j$ this immediately implies that
$$\sum_{i=1}^{\infty}\c_i^{(4)}(x)=0.$$

Finally we have as before by using (\ref{ci5*})
\bea\sum_{i=1}^{\infty}\c_i^{(5)}(x)&=&-\frac{1}{(\a+1)(\a+3)}
\sum_{k=0}^{\a+2}(-1)^k\frac{(-\a-2)_k(\a+3)_k}{(2)_k(\a+4)_kk!}\times{}\nn
& &{}\hspace{1cm}{}\times\sum_{i=0}^{\infty}\sum_{j=0}^i
\frac{(-1)^{i+j}}{\G(i-k)}{k+\a+3 \ch j}{k+\a+3 \ch i-j}(\a+3)_{i-j}x^j.
\n\eea
By using (\ref{K}) and (\ref{Kspec4}) we find
\bea & &\sum_{i=0}^{\infty}\sum_{j=0}^i
\frac{(-1)^{i+j}}{\G(i-k)}{k+\a+3 \ch j}{k+\a+3 \ch i-j}(\a+3)_{i-j}x^j\nn
&=&K(k+3,k+3,3;k;\a,x)=(-1)^{\a+k+1}(\a+3)_{k+1},\;k=0,1,2,\ldots.\n\eea
Hence by using (\ref{som}) again
$$\sum_{i=1}^{\infty}\c_i^{(5)}(x)
=\frac{(-1)^{\a}}{\a+1}\hyphyp{2}{1}{-\a-2,\a+3}{2}{1}
=\frac{(-1)^{\a}}{\a+1}\frac{(-\a-1)_{\a+2}}{\G(\a+4)}=0.$$

From (\ref{ci0}) and the fact that $\c_i(x)=0$ for $i>4\a+10$ we
conclude that
$$\sum_{i=1}^{\infty}\c_i(x)=\sum_{i=1}^{4\a+10}\c_i(x)=0$$
for all $\a\in\{0,1,2,\ldots\}$.

\section{Remarks}

In this section we list some facts we encountered during the research.
First of all we remark that
$$L_n^{(-n)}(-x)=\frac{x^n}{n!}=\sum_{k=0}^n(-1)^k{n+\a \ch n-k}
L_k^{(\a)}(x),\;\ndots.$$
This formula can also be used for inversion instead of lemma~5 in the
following way (compare with lemma~5)~: Suppose that for
$m\in\{0,1,2,\ldots\}$ we have the system of equations
$$\sum_{i=1}^{\infty}A_i(x)D^{i+m}\L=F_n(x),\;n=m+1,m+2,m+3,\ldots,$$
where $\l\{A_i(x)\r\}_{i=1}^{\infty}$ are independent of $n$. Then we
simply find that
$$\sum_{i=1}^{\infty}A_i(x)D^{i+m}\frac{x^n}{n!}=
\sum_{k=m+1}^n(-1)^k{n+\a \ch n-k}F_k(x),\;n=m+1,m+2,m+3,\ldots.$$
Since the generating function
$$e^{xt}=\sum_{n=0}^{\infty}\frac{x^n}{n!}t^n$$
has the trivial inverse $e^{-xt}$ we now conclude that
$$A_i(x)=(-1)^m\sum_{j=1}^i\frac{(-x)^{i-j}}{(i-j)!}
\sum_{k=1}^j(-1)^k{m+j+\a \ch j-k}F_{m+k}(x),\;i=1,2,3,\ldots.$$
However, for our purposes this method turned out to be very inconvenient.

Further we have from (\ref{AAA})
\bea A_0&=&1+M{n+\a \ch n-1}+
\frac{n(\a+2)-(\a+1)}{(\a+1)(\a+3)}N{n+\a \ch n-2}+{}\nn
& &{}\hspace{5cm}{}+\frac{MN}{(\a+1)(\a+2)}{n+\a \ch n-1}{n+\a+1 \ch n-2}.
\n\eea
The coefficients $a_0(n,\a)$, $\b_0(n,\a)$ and $\c_0(n,\a)$ are
respectively
$$a_0(n,\a)=\sum_{k=0}^n{k+\a \ch k-1}={n+\a+1 \ch n-1},\;\ndots,$$
$$\b_0(n,\a)=\sum_{k=1}^n\frac{k(\a+2)-(\a+1)}{(\a+1)(\a+3)}{k+\a \ch k-2}
=\frac{n(\a+2)-\a}{(\a+1)(\a+4)}{n+\a+1 \ch n-2},\;n=1,2,3,\ldots$$
and
$$\c_0(n,\a)=\frac{1}{(\a+1)(\a+2)}
\sum_{k=1}^n{k+\a \ch k-1}{k+\a+1 \ch k-2},\;n=1,2,3,\ldots.$$

Note that the equations (\ref{a1}) and (\ref{a2}) can be written in the form
$$(-1)^k\sum_{i=0}^{\infty}a_i(x)D^{i+k}\L={n+\a \ch n-k}D^2\Ln,\;k=0,1.$$
However it is even more remarkable that the equations (\ref{c1}), (\ref{c2})
and (\ref{c3}) can be written in the form
\bea (-1)^{k+1}\sum_{i=0}^{\infty}c_i(x)D^{i+k}\L&=&
\frac{n(\a+k)-(\a+1)}{(\a+1)(\a+k+1)}{n+\a \ch n-2}
{n+\a \ch n-k}D^2\Ln+{}\nn
& &{}\hspace{1cm}{}+\frac{2}{\a+1}{n+\a \ch n-2}{n+\a \ch n-k}D^3\Ln+{}\nn
& &{}\hspace{1cm}{}+\frac{1}{\a+1}{n+\a \ch n-k-1}
\sum_{i=0}^{\infty}a_i(x)D^{i+2}\L+{}\nn
& &{}\hspace{1cm}{}+{n+\a \ch n-k}\sum_{i=0}^{\infty}b_i(x)D^{i+1}\L,
\;k=0,1,2.\n\eea

Further we remark that the inversion formula of lemma~5 can also be applied
to (\ref{aa1}) instead of (\ref{aa2}) to obtain, by using (\ref{F}) and
(\ref{G})
$$a_i(x)=(-1)^i\l[(\a+1)F_{i-1}(\a+2,\a+2;2;x)-
\frac{1}{\a+2}G_i(\a+2,\a;0;x)\r],\;i=1,2,3,\ldots,$$
which can be shown to be equal to (\ref{ai}) by straightforward but tedious
computation.

Finally we remark that it can be shown that for $\a>-1$ and
$\ell\in\{0,1,2,\ldots\}$
\bea & &\sum_{i=1}^{\infty}a_i(x)D^{i+\ell}\L\nn
&=&-{n+\a \ch n}\sum_{k=0}^{\infty}
\frac{(-n)_{k+\ell+1}(\a+3)_k}{(2)_k(\a+1)_{k+\ell+1}k!}
\hyp{3}{2}{-n+k+\ell+1,-\a-2,\a+k+3}{k+2,\a+k+\ell+2}{1}x^{k+1},\n\eea
for $\ndots$.

\end{document}